\documentclass[11pt]{article}

\usepackage{amsfonts,calrsfs}

\def\Z{\mathbf{Z}}
\def\Q{\mathbf{Q}}
\def\F{\mathbf{F}}
\def\Fr{\mathcal{F}}
\def\x{\mathbf{x}}
\def\y{\mathbf{y}}

\def\C{{\mathcal C}}
\def\D{{\mathcal D}}
\def\uu{\mathbf{u}}
\def\G{\mathcal{G}}

\def\io{\iota}

\newcommand{\vv}{\vspace{1ex}}

\def\pf{{\indent\textit{Proof.}\ }}
\def\qed{\hfill$\square$}
\newcounter{para}[section]
\renewcommand{\thepara}{\thesection.\arabic{para}}
\renewcommand{\thesection}{\arabic{section}}
\renewcommand{\paragraph}{\refstepcounter{para}
\indent{\bf{\thepara}} \ }
\def\sectioning#1{\vv 

 \refstepcounter{section}
\indent {\bf \thesection. #1}}

\begin{document}

{\footnotesize \noindent  Running head: Theory of $\Q$ and elliptic curves \hfill 
March 28th, 2006 

\noindent Math.\ Subj.\ Class.\ (2000): 03B25, 11U05  }

\vspace{30mm}

\begin{center} 

{\Large Elliptic divisibility sequences and 

\vv

undecidable problems about rational points}

\vv

{\sl by} Gunther Cornelissen {\sl and} Karim Zahidi

\vv

\end{center}

\vv

\vv

\begin{quote}
{\small {\bf Abstract.} 
Julia Robinson has given a first-order definition of the rational integers $\Z$ in the rational numbers $\Q$ by a formula $(\forall \exists \forall \exists)(F=0)$ where the $\forall$-quantifiers run over a total of 8 variables, and where $F$ is a polynomial. This implies that the $\Sigma_5$-theory of $\Q$ is undecidable. We prove that a conjecture about elliptic curves provides an interpretation of $\Z$ in $\Q$ with quantifier complexity $\forall \exists$, involving only one universally quantified variable. This improves the complexity of defining $\Z$ in $\Q$ in two ways, and implies that the $\Sigma_3$-theory, and even the $\Pi_2$-theory, of $\Q$ is undecidable (recall that Hilbert's Tenth Problem for $\Q$ is the question whether the $\Sigma_1$-theory of $\Q$ is undecidable). 

In short, granting the conjecture, there is a one-parameter family of hypersurfaces over $\Q$ for which one cannot decide whether or not they all have a rational point.

The conjecture is related to properties of elliptic divisibility sequences on an elliptic curve and its image under rational 2-descent, namely existence of primitive divisors in suitable residue classes, and we discuss how to prove
weaker-in-density versions of the conjecture and present some heuristics.}
\end{quote}

\vv

\vv


{\bf Introduction.} 

\vv

This paper addresses a mixture of number theory and logic, and we will use this introduction to give an informal preview directed at both communities. The central two questions can be phrased as follows: ``What is more difficult: to decide of an arbitrary polynomial equation with integer coefficients whether it has an integer solution, or whether it has a rational solution?''; and: ``What is a hard problem about rational points?'' If one makes these vague questions mathematically more precise, ``decide'' should mean the existence of an algorithm on a Turing Machine (which in practice is equivalent to any notion of ``computable'' via Church's Thesis). Call Hilbert's Tenth Problem ${\rm HTP}(R)$ for a subring $R$ of the rational number $\Q$ the question whether one can decide if an arbitrary polynomial equation with integer coefficients has a solution in $R$. The classical result of Davis, Matijasevich, Putnam and Robinson (\cite{Davis:73}, \cite{M1}, \cite{M2}) shows that ${\rm HTP}(\Z)$, for $\Z$ the ring of integers, has a negative answer. The answer to ${\rm HTP}(\Q)$, however, is not known. But a more general problem has been settled by Julia Robinson in 1949 (\cite{Robinson:49}). She showed that $\Z$ is definable 
in $\Q$ by a first-order formula. This implies that the full first order theory of $\Q$  is undecidable, i.e., that one cannot decide (in the above sense) the truth of an arbitrary first-order sentence in $\Q$ built from the symbols $(0,1,+,\times,=)$. One should think of such a sentence as a ``algorithmically hard'' number theoretical statement 
$$ (\forall x^{(1)}_{1}  \dots x^{(1)}_{f_1})( \exists y^{(1)}_{1} \dots y^{(1)}_{e_1})   \cdots (\forall x^{(N)}_{1}   \dots x^{(N)}_{f_N})( \exists y^{(N)}_{1}   \dots y^{(N)}_{e_N}) \ : \ F({\bf x}, {\bf y}) = 0, $$
where $F$ is a polynomial over $\Z$ in multi-variables ${\bf x}= (x_{1}^{(1)},\dots,x^{(N)}_{f_N})$ and
${\bf y}= (y_{1}^{(1)},\dots,y_{e_N}^{(N)})$.
Note: any formula over $\Q$ can be put into this form, which we call {\sl positive prenex} form (cf.\ lemma \ref{prenex}). Examples of such statements: if there are only existential quantifiers ($N=1,f_1=0$), such a formula says that a certain diophantine equation has a solution; a formula with $N=1$ says that a family of diophantine equations has a solution in ${\bf y}$ for all values of the parameters ${\bf x}$, etc. 

Related to our first question, Robinson's result expresses in some sense that testing the truth of such sentences in $\Q$ or in $\Z$ is ``equally hard''. ${\rm HTP}(\Q)$ is the particular case where one only wants to decide the truth of formul{\ae} with $N=1$ and $f_1=0$ (with $e_1=m$ arbitrary):
$ (\exists y^{}_1 \dots y^{}_{m}) \ : \ F(y_1,\dots,y_m) = 0. $
We now recast the original question above in the following way: how ``complex'' does a formula in $\Q$ have to be, in order to be undecidable? Phrased more dramatically: what is the easiest {\sl hard} problem about rational points? Since we want to indicate how far a formula is from being ``diophantine'' (i.e., in positive prenex form with $N=1$), in \ref{meas}---\ref{compl-reduce} we look at the following two measures of complexity. 
First, we define the {\sl positive} arithmetical hierarchy $(\Sigma^+,\Pi^+)$ as follows: we let $\Sigma^+_0=\Pi^+_0$ denote the set of {\sl atomic} formul{\ae} (=``polynomials''). Define a formula $\mathcal F$ inductively to be in $\Sigma^+_n$ (resp.\ $\Pi^+_n$) if it is of the form $\exists \mathcal G$ (resp.\ $\forall \mathcal G$) with $\mathcal G \in \Pi^+_{n-1}$ (resp.\ $\mathcal G \in \Sigma^+_{n-1}$). The place in the hierarchy of a formula counts its number of {\sl quantifier changes}. Secondly, we introduce the {\sl total number of universal quantifiers} of a formula as above to be $f_1+\cdots+f_N$.

An analysis shows that a positive prenex form of Julia Robinson's original formula defining $\Z$ in $\Q$ is a $\Pi_4^+$-formula (see \ref{compl-julia-nf}), and we can conclude from this that the $\Sigma_5^+$-theory of $\Q$ ($=$ theory of all $\Sigma_5^+$-sentences that are true in $\Q$) is undecidable. ${\rm HTP}(\Q)$ is the question whether the $\Sigma_1^+$-theory is undecidable. Also, that formula, in positive prenex form, has 8 universal quantifiers. This should be considered at the verge of human mathematical understanding --- one is compelled to quote Hartley Rogers, Jr.: ``The human mind seems limited in its ability to understand and visualize beyond four or five alternations of quantifier.'' (\cite{Rogers}, p.\ 322).  

\vv

So how can we, in any way, improve upon this complexity? We propose to use elliptic curves and give a conjectural improvement. First of all, we recall the concept of a {\sl model}\footnote{The words ``model'' and ``interpretation'' seem to have acquired a non-standard meaning in connection with HTP. The precise meaning will be explained in the text.}  of $\Z$ in $\Q$ (cf.\ \ref{diomodel}) and study how the complexity of formul{\ae}  changes under interpretation via certain models (\ref{def-compldiomodel}--\ref{undecQ}). 
We then recall (in Section \ref{eds}) how elliptic curves over $\Q$ provide natural models of $(\Z,+)$ in $\Q$. We follow a suggestion of Pheidas (\cite{Pheidas}) that it is natural to use such models to try to define  ``divisibility'' of integers within $\Q$; this is very much inspired by the function field case. For this, we have to introduce a variant of the old concept of ``elliptic divisibility sequence'' (apparently due to Lucas and studied by M.~Ward, cf.\ \cite{Ward:48}). Assume that $E$ is an elliptic curve over $\Q$ with $(0,0)$ as 2-torsion point and Weierstrass equation $y^2=x^3+ax^2+bx$ with $b$ squarefree, and that $P$ is a point of infinite order of sufficiently large height on $E$. Then for even $n$, we can write $$nP=(x_n,y_n)=\left( \left( \frac{A_n}{B_n}\right)^2,\frac{A_n C_n}{B_n^3}\right)$$ for coprime integers $A_n,B_n,C_n$, and $\{C_\ast\}$ forms an {\sl odd divisibility sequence} in the sense that $C_n$ divides $C_{tn}$ precisely for odd $t$ (\ref{odd}). 
Our first main theorem uses two further notions. Let $R$ denote a set of primes. We agree to identify primes $p$ with normalized non-archimedean valuations $v=v_p$, such that $v(p)=1$ and $v(ab)=v(a)+v(b)$ (please mind: this is a logarithm of what has been called a valuation elsewhere). We say that $\{ C_\ast \}$ is $R$-(odd-)primitive if any $C_n$ has an (odd order) primitive divisor from $R$, i.e., there is a valuation $v \in R$ such that $v(C_n)$ is non-zero (odd) but $v(C_i)=0$ for all $i<n$. Secondly, for two rational numbers $x$ and $y$, we denote by ${\mathcal D}_R$ the {\sl $R$-divisibility predicate}: $(\forall v \in R)(v(x) \mbox{ odd } \Rightarrow v(x) < v(y^2)$). Theorem \ref{defdiv} then says that {\sl if in the above setup, $\{ C_\ast \}$ is $R$-odd-primitive, then for any integers $m,n \in {\bf Z}$,
\begin{eqnarray*} 
m | n  &\iff& \D_R(y_m \sqrt{x_m}, y_n \sqrt{x_n}) \vee \D_R(y_m \sqrt{x_m}, y_{m+n} \sqrt{x_{m+n}})
\end{eqnarray*}}
This is our attempt at defining integer divisibility in the rational numbers.

The relevant question becomes: can we find $R$ for which ${\mathcal D}_R$ is equivalent to a formula in $\Sigma_1^+$ (whence irrelevant from our point of view of complexity) and for which $\{ C_\ast \}$ is $R$-primitive? The elliptic Zsigmondy's theorem, transferred to $C$, says that $R$-primitivity holds for $R$ equal to the set of all primes, but we don't know whether ${\mathcal D}_R$ is $\Sigma_1^+$ for that $R$. On the other hand, a theorem of Van Geel and Demeyer (based on previous work of Pheidas and Van Geel/Zahidi) states that ${\mathcal D}_R$ is diophantine for $R=R_D$ the set of primes inert in one of finitely many quadratic number fields of discriminants $D=\{d_1,\dots,d_r\}$, and hence for $R$ of arbitrary high Dirichlet density $\neq 1$, see \ref{dir}. So our natural conjecture (\ref{SC}) becomes an {\sl inertial elliptic Zsigmondy's theorem: there exists $E,P$ and $D$ as above such that $\{C_\ast\}$ is $R_D$-odd-primitive.}

We can show that multiplication is definable in $(\Z,+,|,0,\neq)$ by a $\Sigma_3^+$-formula only involving one universal quantifier (\ref{lemdiv}), and that our model allows us to get rid of ``$0$'' and ``$\neq$''.
Collecting these facts, we arrive at our second main theorem: {\sl the conjecture implies that integer arithmetic $(\Z,+,\times)$ is interpretable by a $\Sigma_3^+$-formula in the rationals $\Q$, using only one universal quantifier; and that the $\Sigma_3^+$-theory of $\Q$ is undecidable.} This (conjecturally) improves the complexity of Robinson's definition in two ways. In section \ref{pi2}, we adapt the construction to show that {\sl the conjecture even implies that the $\Pi_2^+$-theory (and even the set of formul{\ae} with only one universal quantifier) is undecidable}; but note that this is {\sl not} proven by constructing a model of $\Z$ in $\Q$ that has complexity $\Pi_2^+$. The geometrical meaning of this statement is that  there is a one-parameter family of hypersurfaces over $\Q$ for which one cannot decide whether or not they all have a rational point.

It is difficult to verify the conjecture numerically since it involves hard prime factorisations. However, note that the philosophy of encoding the integer $n$ by the point $nP$ on an elliptic curve is advantageous from the point of view of divisibility for two reasons: the ``powerful'' part of the coordinates of $nP$ is very small (in the sense that the height of the ``powerless'' part of $nP$ is of the same order as the height of $nP$), and $C_n$ tends to have many more prime factors than $n$. These remarks can be turned into heuristics that support the conjecture (see section \ref{conj}). The conjecture incorporates statements about solutions in coprime integers of such Calabi-Yau surfaces as $$ (A^2+B^2)(A^2+11B^2)=3^2 \cdot 5^2 \cdot (X^2-5Y^2)^2, $$
which becomes the ``One Equation to Rule Them All'' of the $\Pi_2$-theory of $\Q$ (like Martin Davis's for the $\Sigma_1$-theory of $\Z$; but that equation heuristically behaved the wrong way, cf.\ \cite{Shanks2}).

Finally, we use the periodicity of elliptic divisibility sequences to prove in Section \ref{densversion} that {\sl if $\{ B_\ast \}$ is the elliptic divisibility sequence associated to $(2,-4)$ on $y^2=x^3+7x^2+2x$, then any $B_{s^e}$ for $s$ a prime number $=\pm 3$ mod 8 has a {\sl primitive} odd order divisor from $R_5$, and for $D=\{5,13,29,41,53\}$, the set $\{ s \mbox{ prime} : B_s \mbox{ has a primitive odd order divisor from } R_D \}$ has Dirichlet density at least $95.5 \% $.}

\vv

{\bf Remarks.} (i) Beltjukov studied the theory of $(\Z,+,|)$ (\cite{Bel}) and Lipshitz (\cite{Lip2}, \cite{Lip}, \cite{Lip1}) has studied divisibility structures of the form $({\mathcal O},+,|)$ for $\mathcal{O}$ the ring of integers in a number field $K$, including (independently) the usual integers, and obtained exact results on which of these theories are (un)decidable. He showed in particular that multiplication is definable in the $\Sigma^+_1$-theory of such a structure and that it contains a diophantine model of $\Z$, precisely if ${\mathcal O}$ has infinitely many units. Thus, if $K$ is not equal to $\Q$ or an imaginary quadratic number field and $A$ is an abelian variety with multiplication by ${\mathcal O}$, then an imitation of the above theory for $A$ would lead to a $\Sigma_1^+$-definition of $\Z$ in $\Q$ and hence a negative answer to Hilbert's Tenth Problem for $\Q$. 
This can already occur for $A$ the Jacobian of a genus two curve with real multiplication: to give an example from \cite{Stoll}, the curve $$ \mathcal{C} \ : \ y^2 + (x^3+x^2+x)y=x^4+x^3+3x^2-2x+1 $$
(the modular curve $X_0(85)$ modulo an Atkin-Lehner involution)
has a Jacobian of rank two over $\Q$, and real multiplication by $\Z[\sqrt{2}]$ defined over $\Q$.  The r\^ole of the ``$x$-coordinate'' on the elliptic curve should be played by the associated Kummer surface.
There are, however, many obstacles to make such a generalisation work, even assuming certain arithmetical 
conjectures. 

A generalisation of the above to elliptic curves with complex multiplication, however, should be unproblematic.

\vv

(ii) In another direction, Poonen (\cite{Poonen:2003}) has shown that there exists a set $S$ of primes of Dirichlet density one such that $\Z$ is definable by a diophantine formula in $\Z[\frac{1}{S}]$. 

\vv

(iii) This paper supersedes the first author's year 2000 manuscripts \cite{Cornelisseneds} about the topic.

\vv

(iv) Number theorists can take the following direct path to the relevant conjecture: Sections \ref{divpoly}-\ref{exa} (divisibility sequences), \ref{primcond}, \ref{weakprim} ((weak) $R$-primitivity), \ref{SC} (main conjecture) and Sections \ref{conj} and \ref{densversion} (discussion of the conjecture). 

\vv

\vv

{\small {\bf Acknowledgements.} \ It is a pleasure to thank Thanases Pheidas for his help. The first author thanks Graham Everest for making him reconsider this material during an inspiring visit to UEA in 2004, and for many suggestions. Some of the heuristical arguments in section \ref{heur} were shown to us by Bjorn Poonen at Oberwolfach in 2003, and some of the references in that section where kindly provided by Pieter Moree. Marco Streng suggested some important improvements in section \ref{eds}.}

\vv

\vv

\sectioning{Models and their complexity.} \label{compl} 

\vv

\paragraph {\bf Positive prenex-form.} \label{compl-intro} \ Julia Robinson proved in 1949 that the set of rational integers $\Z$ is definable in the rational numbers $\Q$ (\cite{Robinson:49}) by a first-order formula. It is still
an open problem whether $\Z$ can be defined in $\Q$ by a positive-existential formula (and consequently, the positive-existential theory of $\Q$ is undecidable). It should therefore be interesting to study the question of how complicated 
the definition of $\Z$ is in terms of number of universal quantifiers used, 
or number of quantifier changes. We thus propose to study the complexity of defining the integers in the rational numbers. To formulate the problem very
precisely, we need to make the following convention: a formula $\Fr$ in the
first-order theory of $(\Z,+,\times,0,1,=)$ or $(\Q,+,\times,0,1,=)$ will be written in the following normal
form:
$$ \forall x^{(1)}_1 \dots \forall x^{(1)}_{f_1} \exists y^{(1)}_1 \dots \exists y^{(1)}_{e_1} \cdots \forall x^{(N)}_1 \dots \forall x^{(N)}_{f_N} \exists y^{(N)}_1 \dots \exists y^{(N)}_{e_N} \ : \ F({\bf x}, {\bf y}) = 0, $$
with $e_i>0$ for $i=1,\dots,N-1$ and $f_i>0$ for all $i=2,\dots,N$; where $F$ is a polynomial in multi-variables ${\bf x}= (x^{(1)}_1,\dots,x^{(1)}_{f_1},\cdots,x^{(N)}_1,\dots,x^{(N)}_{f_N})$ and
${\bf y}= (y^{(1)}_1,\dots,y^{(1)}_{e_1},\cdots,y^{(N)}_1,\dots,y^{(N)}_{e_N})$.
We will call such a formula {\sl a $((f_1,e_1),\dots,(f_N,e_N))$-formula} and call this form the {\sl positive prenex} form. Note that the formula is not only in ``prenex''-form (in which the quantifiers are followed by a quantifier-free formula that can be any boolean combination of atomic formul{\ae}; see, e.g. \cite{CL}, p.\ 157), but that we let the quantifiers be followed by a single atomic formula, viz.\, an equation. That this is possible is 
specific to certain languages. We don't want to allow negations in the quantifier-free part, because we are interested in measuring ``closeness'' 
to a {\sl positive existential} ($=$ diophantine) formula. It is indeed possible to transform any formula into such positive prenex normal form; this is well known but we include a proof for completeness.

\vv

\paragraph {\bf Lemma.} \label{prenex} {\sl let $R \subseteq \Q$ be a ring. Any first-order formula in the ring language $(R,+,\times,0,1,=)$ of $R$ can be written in normal form.}

\vv

\pf The following logical connectives can occur: $ \Rightarrow, \neg, \vee, \wedge$. Here is an algorithm that eliminates their occurrences. Replace $A \Rightarrow B$ by $\neg A \vee B$. Pull negations from left to right through a formula (changing quantifiers and connectives accordingly). Put all the quantifiers on the left (possibly changing names of variables). 

Lagrange's four-squares theorem states that any integer $n \geq 0$ is
a sum of four squares. Therefore, for $x \in R$ we have $$x>0 \iff (\exists a,b,c,d,e,f,g,h)((e^2+f^2+g^2+h^2)(x+1)=a^2+b^2+c^2+d^2).$$ Furthermore, $n \neq 0 \iff (n>0)\vee(n<0)$. Use this to replace, for a polynomial $P$, the formula $P \neq 0$ by a formula only involving equality signs. 

For polynomials $P$ and $Q$, replace $(P=0) \vee (Q=0)$ by $PQ=0$, and  $(P=0) \wedge (Q=0)$ by $P^2+Q^2 = 0$. The final result of all these replacements is the above normal form. \qed

\vv

\paragraph {\bf Remark.} \ \label{compl-rmkQ} Depending on $R$, one can sometimes improve upon the number of existential quantifiers used to translate $P \neq 0$. For example, if $R=\Q$, then $P \neq 0 \iff (\exists Q)(PQ=1)$. 

\vv

\paragraph {\bf Remark.} \ The polynomial $F$ in the general positive prenex form might still depend on unquantified variables (also called
free variables) which are omitted in our notation; this will cause no
confusion. If no free variables occur we call the formula a {\sl sentence}. A
sentence has a precise truth-value, whereas this is not the case for a
formula with free variables. However if we give these free variables a
specific value then we obtain a sentence with a specific truth value. The
set of all specifications of the free variables for which the
corresponding sentence is true, is {\sl the set defined by the formula}.

\vv

\paragraph {\bf Measures of complexity.} \ \label{meas} As explained in the introduction, we do not care too much about the number of existential quantifiers, but want to have as few universal quantifiers as possible in our formul{\ae}. A first measure of such complexity of a formula is its {\sl total number of universal quantifiers} ($t$-complexity) $$t(\Fr):=f_1 + \cdots + f_N. $$ A second measure of complexity is the place of the formula in the (positive) arithmetical hierarchy, that we will now introduce. 

\vv

\paragraph {\bf The positive arithmetical hierarchy.} \label{defpisigma} One usually defines the (arithmetical) hierarchy $(\Sigma, \Pi)$ of a language as follows (compare \cite{CK}, p.\ 117). Let $\Sigma_0=\Pi_0$ denote the set of quantifier-free formul{\ae}. Define a formula $\mathcal F$ inductively to be in $\Sigma_n$ (resp.\ $\Pi_n$) if it is of the form $\exists \mathcal G$ (resp.\ $\forall \mathcal G$) with $\mathcal G \in \Pi_{n-1}$ (resp.\ $\mathcal G \in \Sigma_{n-1}$). 

In accordance with our use of a normal form which is positive prenex, we define the {\sl positive arithmetical hierarchy} $(\Sigma^+,\Pi^+)$ as follows: we 
let $\Sigma^+_0=\Pi^+_0$ denote the set of positive boolean combinations of atomic formul{\ae}. Define a formula $\mathcal F$ inductively to be in $\Sigma^+_n$ (resp.\ $\Pi^+_n$) if it is of the form $\exists \mathcal G$ (resp.\ $\forall \mathcal G$) with $\mathcal G \in \Pi^+_{n-1}$ (resp.\ $\mathcal G \in \Sigma^+_{n-1}$). 

A formula in $\Sigma_1$ is called {\sl existential}, in $\Pi_1$ {\sl universal}, in $\Sigma^+_1$ {\sl positive existential} or {\sl diophantine}.

The {\sl number of quantifier changes} $c$ ($c$-complexity) can be defined by $$ c(\Fr) := \left\{ \begin{array}{ll}  2N-1 & \mbox{if } f_1e_N \neq 0, \\ 2N-2 & \mbox{if one of } f_1, e_N = 0 \\ 2N-3 & \mbox{if } f_1=e_N=0. \end{array} \right. $$    
In terms of the hierarchy, this means the following: if ${\mathcal F} \in \Sigma^+_{n+1}-\Pi^+_{n}$ or ${\mathcal F} \in \Pi^+_{n+1}-\Sigma^+_{n}$, then $c({\mathcal F})=n$. 

For a ring language as in \ref{prenex}, formul{\ae} in $\Sigma_0^+$ are equivalent to {\sl atomic} formul{\ae} by \ref{prenex}. Furthermore, as non-equalities are existential, any $\Sigma_{2n+1}$-formula is equivalent to  a $\Sigma_{2n+1}^+$-formula and any $\Pi_{2n}$-formula is equivalent to a $\Pi_{2n}^+$-formula. 

By abuse of the syntax/semantics difference, we will from now on sometimes write that $\Fr \in \Sigma_n^+$ if $\Fr$ is equivalent in the theory under consideration to a formula in  $\Sigma_n^+$.

\vv

\paragraph {\bf Remark.} \label{compl-davis} \ (i) For $({\bf N},+,\times,0,1)$, a polynomial bijection ${\bf N}^2 \rightarrow {\bf N}$ as in Martin Davis (\cite{Davis:73}, pp. 236-237) can be used to show that any formula is equivalent to a formula in positive prenex form with $f_1=\dots=f_N=1$. For  $({\bf Z},+,\times,0,1)$, the same conclusion $f_1=\dots=f_N=1$ holds by the method of diophantine storing. The analogous statement is not known for $\Q$, but would follow from the ABC-hypothesis, see \cite{Cor}. 

\vv

(ii) In the course of the proof of the main theorem, we will also have to use other languages than the usual ring language, and the reader should be cautioned that the positive and the usual hierarchy can be quite different in such a case (up to equivalence of formul{\ae} in that language): there might be quantifier-free formul{\ae} that are not equivalent to an atomic formula. Example: $(a=b) \wedge (c=d)$ in $(\Z,+,0,1,=)$.

\vv

\paragraph {\bf Remark.} \label{compl-reduce} \ A
formula $\Fr$ could be equivalent (in a given theory) to a formula $\mathcal G$ whose complexity is different. In practice, it is often possible to reduce the number of universal quantifiers in a formula by using fewer variables, and we will sometimes do so. For example if $\Fr$ and $\mathcal G$ are formul{\ae} with disjoint sets of variables, then $(\forall X,Y)(\Fr(X) \wedge {\mathcal G}(Y))$ is equivalent to $(\forall X)(\Fr(X) \wedge {\mathcal G}(X))$. 

\vv

\paragraph {\bf Robinson's definition.} \ Julia Robinson's definition of the integers is the following: Let $\phi(A,B,K)$ denote the formula $(\exists X,Y,Z)(P_{A,B,K}^{X,Y,Z}=0)$ with $P_{A,B,K}^{X,Y,Z} = 2+ABK^2+BZ^2 -X^2 - AY^2$. Then for $N \in \Q$, we have $N \in \Z  \iff \mathcal{R}(N)$ with 
\begin{eqnarray*}
\mathcal{R}(N) & : &    \forall A,B  \{ \ [ \ \phi(A,B,0) \ \wedge \ (\forall M)(\phi(A,B,M) \Rightarrow \phi(A,B,M+1)) \ ] \\
& & \Rightarrow \phi(A,B,N)  \ \}
\end{eqnarray*}

\noindent We will now analyse the diophantine complexity of this formula:

\vv

\paragraph {\bf Lemma.} \label{compl-julia-nf} {\sl The formula $\mathcal{R}$ is equivalent to a  $\Pi_4^+$- $((5,4),(3,1))$-formula $\Fr$ with $t(\Fr)=8$ and $c(\Fr)=3$.}

\vv

\pf We use the algorithm  from the proof of Lemma \ref{prenex}. Thus, we replace the implications
to get
$$(\forall A,B) \{ \  \neg [ \ \phi(A,B,0) \ \wedge \ (\forall M)(\neg \phi(A,B,M) \vee \phi(A,B,M+1)) \ ]  \vee \phi(A,B,N) \ \} .$$
We pull through the negations
$$(\forall A,B) \{ \   [ \ \neg \phi(A,B,0) \ \vee \ (\exists M)( \phi(A,B,M) \wedge \neg \phi(A,B,M+1)) \ ]  \vee \phi(A,B,N) \ \} .$$
Now plug in the $(0,3)$-formula $\phi(A,B,\ast)$  
\begin{eqnarray*} & & (\forall A,B,X,Y,Z)( \exists M,X',Y',Z')( \forall X'',Y'',Z'' )(\exists X''',Y''',Z''') \\ & & [ \ P_{A,B,0}^{X,Y,Z} \neq 0 \vee (P_{A,B,M}^{X',Y',Z'}=0 \wedge P_{A,B,M+1}^{X'',Y'',Z''} \neq 0) \vee P_{A,B,N}^{X''',Y''',Z'''}=0 \ ]
\end{eqnarray*}
In $\Q$, we can replace an inequality by an equality at the cost of introducing one existential quantifier (\ref{compl-rmkQ}). We can use the same variable for both inequalities in the above formula, since it is a disjunction of inequalities. We can simplify the arising formula further by using the same name for $X'$ and $X'''$, $Y'$ and $Y'''$ and $Z'$ and $Z'''$. to arrive at a $((5,4),(3,1))$-formula. \qed

\vv

\paragraph {\bf (Diophantine) models.} \label{diomodel} \ 
Our (conjectural) improvement of this formula will not depend on a definition of $\Z$ as a {\sl subset} of $\Q$, but rather on the existence of a model of $\Z$ over $\Q$. We therefore give a general definition first (in a certain model theoretic parlance, this just means an interpretation of the first theory in the second model):

\vv

\paragraph {\bf Definition.} \label{def-diomodel} \ Let $(M,L,\phi)$ be a triple consisting of a set $M$ and a finite collection $L=\{r_i\}$ of subsets of cartesian powers of $M$ (called ``relations'' or ``constants''), where $\phi$ is an {\sl interpretation} of $L$ in $M$ (which we will often leave out of the notation). If $(N,L'= \{ s_i \},\phi')$ is another such triple, $M$ is said to have a {\sl model $(D,\iota)$ in $N$} if there is a bijection $\iota  :  M \rightarrow D$ between $M$ and a definable subset $D$ of some cartesian power $N^d$
of $N$, such that the induced inclusions of $\iota(r_i)$ in the appropriate cartesian power of $N$ are definable subsets. We call $d$ the {\sl dimension} of the model. By slight abuse, we will sometimes omit $\iota$ from notations.

\vv

\paragraph {\bf Examples.}\  \label{ex1-diomodel} From now on, we will write $\Z$ and $\Q$ for $(\Z,L)$ and $(\Q,L)$ with $L=(0,1,+,\times,=)$ the standard language of rings. By further abuse of notation, we will often leave out the constants ``$0$'',``$1$'' and equality ``$=$'' from a language on a ring. A model of $\Z$ in $\Q$ is a countable definable subset of $D$, such that under a bijection $ \iota   :  \Z \rightarrow D$, the induced images of the graphs of addition and multiplication are definable subsets $D_+$ and $D_\times$ of $\Q^3$. The result of Julia Robinson shows that one can cake $D={\bf Z}$ and $\iota=\mbox{id}$, leading to a one-dimensional model. If $G$ is an affine algebraic group over $\Q$, then embedding $G$ in some affine space of dimension $d$ gives a $d$-dimensional model of $(G(\Q),+_G)$ in $(\Q,+,\times)$.  If $G(\Q) = \Z$, one thus has a model of $(\Z,+)$ in $\Q$ (but lacking multiplication).  

\vv

One can measure the complexity of a model by the complexity of the formul{\ae} that define the embeddings of the relations. Thus, 

\vv

\paragraph {\bf Definition.} \ \label{def-compldiomodel} For $S$ a definable subset of a cartesian power of $N$, write $t(S) \leq n$ (or $c(S) \leq n$) if there exists a formula $\Fr$ defining $S$ with $n=t(\Fr)$ (or $n=c(\Fr)$). 

We say that the {\sl $t$-complexity $t(D)$ of 
a model} $(D,\iota)$ of $(M,L)$ in $(N,L')$ satisfies $t(D) \leq n$ if $$\max \{ t(\iota(D)), t(\iota(r_i)) \} \leq n,$$
and similarly for the $c$-complexity or position in the hierarchy. $D$ is called a {\sl diophantine model} of $M$ in $N$ if $t(D)=0$. 

\vv

\paragraph {\bf Remark.} \ This definition involves only upper bounds for the complexity of a definable set, since $S$ could be definable by several equivalent formul{\ae} having different complexity, cf.\ \ref{compl-reduce}. In general, it seems quite hard to prove that a set {\sl cannot} be defined by a less complex formula. 

\vv

\paragraph {\bf Examples} (continued). \ \label{ex2-diomodel} The complexity of Julia Robinson's model is as in Lemma \ref{compl-julia-nf}. The $t$-complexity 
of embedding $(G(\Q),+_G)$ in $\Q$ (for $G$ an affine algebraic group) is zero, since $G(\Q)$ is the solution set to the ideal of equations that defines $G$ in affine space, and addition is defined by an algebraic formula that involves the coordinates in that affine space (note that a different formula might be needed for distinct cases, such as doubling of points, but this distinction is made by a formula only involving inequalities and case distinctions, that are equivalent to a formula only involving existential quantifiers). 

\vv

\paragraph {\bf Remark.} \ If $\Z$ admits a diophantine model in $\Q$, then there exists a variety $V$ over $\Q$ such that the real topological closure of the set of rational points $V(\Q)$ in the set of real points $V({\bf R})$ has infinitely many connected components. This contradicts a conjecture of Mazur, cf.\ \cite{CZ}.

\vv

\paragraph {\bf Translation of formul{\ae}.} \ One can use a model of $M$ in $N$ to translate formul{\ae} in $M$ to formul{\ae} in $N$, such that true sentences in $M$ are precisely translated into true sentences in $N$. Given a formula $\Fr$ in $M$, one replaces every occurrence of a variable $x$ by the $N$-definition of ``$x \in D$'', and every occurrence of a relation $r({\bf x})$ by the $N$-definition of $r$.
One thus gets a formula which we denote by $\iota({\Fr})$. 

\vv

\paragraph {\bf Example.} \label{ex-trans} Consider the formula $\Fr \ : \ (\exists x_1)(\forall x_2)(x_1^2 x_2 + x_2 = 0)$ in $\Z$. Suppose one is given a 2-dimensional model $D \subseteq \Q^2$ of $\Z$ in $\Q$. Then this formula translates into
\begin{eqnarray*} \iota(\Fr) & : &  (\exists y_1^1 y_1^2)(\forall y_2^1 y_2^2)(\exists u_1 u_2 v_1 v_2)[(y_1^1,y_1^2) \in D \wedge [(y_2^1,y_2^2) \in D  \Rightarrow \\ 
& & [(y_1^1,y_1^2,y_1^1,y_1^2,u_1,u_2) \in \iota(\times) \wedge  (y_2^1,y_2^2,u_1,u_2,v_1,v_2) \in \iota(\times) \wedge \\ & & 
(y_2^1,y_2^2,v_1,v_2,\iota(0)) \in \iota(+)]]]
\end{eqnarray*}
where one should now further replace membership of $D, \iota(+)$ and $\iota(\times)$ by their first-order definitions. Note the introduction of the ``dummy variables'' $u_i,v_i$  to unravel nested occurrences of addition and multiplication. 

\vv

If one applies positive prenex simplification to remove implications and negations, one can keep track of the complexity of the translation. One can ask how the complexity of a formula changes under translation. We will only consider the following case:

\vv

\paragraph {\bf Proposition.} \ \label{compl-change} {\sl Let $(D,\io)$ be a $d$-dimensional model of $\Z$ in $\Q$ and assume that membership of $D$ 
is atomic and  of $\io (+)$ is $\Sigma _{1}^{+}$. 
In the following table, the second and third column list the positive hierarchical status of the formula $\io(\Fr)$ as a function of the status of $\io(\times)$ and $\Fr$ as it is indicated in the first column and top row:

\renewcommand{\arraystretch}{1.2}

\begin{center}
\begin{tabular}{l|l|l}
$\Fr \in$ & $\io(\times) \in \Sigma_{s}^{+} \Rightarrow \io(\Fr) \in$ & $\io(\times) \in \Pi_{s}^{+} \Rightarrow \io(\Fr) \in$\\
\hline
$\Sigma_{2n}^{+}$ & $\Sigma_{2n+s}^{+}$ & $\Sigma_{2n+s+1}^{+}$ \\ $\Sigma_{2n+1}^{+}$ & $\Sigma_{2n+s}^{+}$ & $\Sigma_{2n+s+1}^{+}$ \\
 $\Pi_{2n}^{+}\; (n>0)$ & $\Pi_{2n+s-1}^{+}$ & $\Pi_{2n+s}^{+}$ \\
 $\Pi_{2n+1}^{+}$ & $\Pi_{2n+s+1}^{+}$ & $\Pi_{2n+s+2}^{+}$ \\
\end{tabular}
\end{center}
(note: inclusion of a formula in a class of the hierarchy means that the formula is equivalent to a formula in that class). 
Furthermore, in all cases we have 
$$t(\io (\Fr)) \leq t(\iota(\times)) +dt(\Fr )\; .$$
 
}
\vv

\pf The proof is a matter of non-trivial book-keeping. We use the following notation: let $\Delta _{n}^{+}=\Sigma _{n}^{+}\cup \Pi_{n}^{+}$. 
We need to establish the following fact, that will be used implicitly in the sequel:

\vv

\vv

{\bf \ref{compl-change}.1 Lemma.} \ {\sl Let  $\Fr_{1}$ be a $\Sigma _{n}^{+}$-formula (respectively a $\Pi _{n}^{+}$-formula) and suppose that $\{ \Fr_{2}, ...,\Fr_{q}\}$ is a finite collection of formul{\ae}   such that each $\Fr_{i}$ is either a 
$\Sigma _{n}^{+}$-formula (respectively a $\Pi _{n}^{+}$-formula) or a $\Delta _{m}^+$-formula, for some $m<n$. Then $$\Fr=\Fr_{1}\wedge...\wedge \Fr_{q} \mbox{ and }\Fr '=\Fr_{1}\vee...\vee \Fr_{q}$$ are $\Sigma _{n}^{+}$-formul{\ae}   (respectively $\Pi _{n}^{+}$-formul{\ae}). 
 
 Suppose further that $\Fr_i$ is a $((f_{ik_i},e_{ik_i}),...,(f_{i1},e_{i1}))$-formula; then:
$$t(\Fr )\leq \sum_{j=1}^{k_1} \max\{f_{1j},...,f_{qj}\}\; 
 \mbox{  \ and \ } t(\Fr ') \leq t(\Fr_1)+...+t(F_q) = \sum_{i=1}^q \sum_{j=1}^{k_i} f_{ij} \; . $$
 }

 We prove the result for $\Fr _1$ a $\Pi_n^+$ statement -- the other cases are similar.
Without loss of generality, we may assume that each formula $\Fr _{i}$ is a $\Pi_{n}^+$-formula (indeed, we can add quantifiers whose variables are those variables that do not appear freely in $\Fr_i$; for each new variable $x$ introduced in this way add the 
equation ``$x=x$'').
For $n=0$ the statement is trivial. For $n=1$ the result is also clear, since for any formul{\ae}   $\C$, $\D$, $(\forall x)(\C(x))\wedge (\forall y)(\D(y))$ is equivalent to $(\forall x)(\C(x)\wedge \D(x))$. So, suppose  each of the $\Fr _{i}$ is a $\Pi _{n+1}^+$-formula with $n>0$, i.e., each $\Fr_i$ is of the form $$(\forall x_{1},...,x_{f_{i1}})(\exists y_{1},...,y_{m_{i}})
(\G _{i}(x_{1},...,x_{n},y_{1},...,y_{f}))$$ with $\G_{i}\in \Pi_{n-1}^{+}$. Let $m=m_{1}+...+m_{q}$ and $f=\max_{i}f_{i1}$. Since for any formul{\ae}   $\C$, $\D$, $(\forall x)(\C(x))\wedge (\forall y)(\D(y))$ is equivalent to $(\forall x)(\C(x)\wedge \D(x))$, and $(\exists x)(\C (x))\wedge (\exists y)(\D(y))$ is equivalent to $(\exists x,y)(\C (x)\wedge (\D(y))$, the formula $\Fr $ is  equivalent to
$$ (\forall x_{1},...,x_{f})(\exists y_{1},...,y_{m})(\G _{1}\wedge ...\wedge \G_{q})\; .$$ 
By induction, the formula  $\G= \G_{1}\wedge ...\wedge \G_{q}$ is 
 $\Pi _{n-1}^{+}$, 
hence $\Fr $ is a $\Pi_{n+1}^{+}$-formula. We have
$t(\Fr)=f+t(\G)$, and hence by induction $t(\Fr)=\max_if_{i1}+\sum_{j=2}^{k}\max_{i}f_{ij}$, which proves the result (note that the extra quantifiers which we may have added to make all formul{\ae}   $\Pi _n^+$ do not affect the statement). The statement concerning a disjunction of $\Pi _{n}^{+}$-formul{\ae}   can be proven similarly, by noting that for formul{\ae}   $\C$, $\D$, $(\forall x)(\C(x))\vee (\forall y)(\D(y))$ is equivalent to $(\forall x,y)(\C(x)\vee \D(x))$, and $(\exists x)(\C (x))\vee (\exists y)(\D(y))$ is equivalent to $(\exists x)(\C (x)\vee (\D(x))$. This proves the lemma. \qed

\vv

\vv

The proof of \ref{compl-change} is by induction, jumping down by 2 in the hierarchy (and thus induction starts at the two lowest levels of the hierarchy): 
\vv

(a) Let $\iota(\times) \in \Sigma^+_s$. 
Suppose first that $\Fr\in \Sigma_0^+$, i.e., $$\Fr: \; F(x_{1},...,x_{n})=0$$ for some integral polynomial.
Then there exists a set $\Lambda=\{1,...,\ell\}$, with $\ell \geq n$, subsets $I,J\subset
 \Lambda ^{3}$ and natural numbers $s,r\in \Lambda$ such that the translation $\io (\Fr)$ is of the form:
$$(\x_{1}\in D\wedge ...\wedge\x_{n}\in D)\wedge (\exists \uu_{1},...,\uu_{\ell})$$
$$(\uu_{1}=\x_{1}\wedge...\wedge\uu_{n}=\x_{n}  \bigwedge_{\bar{i}\in I}(\uu _{i_{1}},\uu _{i_{2}},\uu _{i_{3}})\in \io (+) $$
$$ \bigwedge_{\bar{j}\in J}(\uu _{j_{1}},\uu _{j_{2}},\uu _{j_{3}})\in \io (\times)\wedge (\uu _{r},\uu _{r},\io (0))\in \io (+)) \; ,$$
where $\bar{i}=(i_{1},i_{2},i_{3})$ and $\bar{j}=(j_{1},j_{2},j_{3})$ are multi-indices and boldface variables are variables ranging over $\Q^d$.
The conjunction
$$\bigwedge_{\bar{i}\in I}(\uu _{i_{1}},\uu _{i_{2}},\uu _{i_{3}})\in \io (+) 
 \bigwedge_{\bar{j}\in J}(\uu _{j_{1}},\uu _{j_{2}},\uu _{j_{3}})\in \io (\times)\wedge (\uu _{r},\uu _{r},\io (0))\in \io (+)) $$
is a conjunction of $\Sigma _1^+$-formul{\ae}   and $\Sigma _{s}^+$-formul{\ae}, and hence is itself a $\Sigma_{s}^+$-formula. Hence, the formula  $\io (\Fr)$ is a conjunction of a  $\Sigma_{s}^{+}$-formula and a $\Sigma _0^+$formula, hence is a  $\Sigma_{s}^{+}$-formula. Furthermore, $t(\io (\Fr))=t(\io (\times))$.\\

If $\Fr\in\Sigma_{1}^{+}$, then $\io (\Fr) \in \Sigma_{s}^+$. Indeed, write $\Fr= (\exists x_1,...,x_n)\G(x_1,...,x_n)$ for some quantifier-free formula $\G$. The translation $\io (\Fr)$ is then given by:
 $$(\exists \x_1,...,\x_n)(\x_1\in D\wedge ...\wedge \x_n\in D \wedge \io(\G)(\x_1,...,\x_n)),$$ and the result is clear since membership of $D$ is atomic.  

  We proceed by induction. Suppose that $\Fr$ is  $\Sigma_{n+1}^{+}$, i.e. there exists a 
$\Sigma _{n-1}^+$-formula $\G$ such that $$\Fr:\; (\exists x_{1},...,x_{n})(\forall y_{1},...,y_{m})(\G (x_{1},...,x_{n},y_{1},...,y_{m}))\; .$$
The translation $\io (\Fr)$ then becomes:
\begin{eqnarray*} & & (\exists \x_{1},...,\x_{n})(\forall \y_{1},...,\y_{m})  [(\x_{1}\in D\wedge ...\wedge\x_{n}\in D)\wedge  \\ & &  \wedge ((\y_{m}\in D\wedge ...\wedge\y_{m}\in D)
\Rightarrow \io (\G)(\x_{1},...,\x_{n},\y_{1},...,\y_{m}))] \end{eqnarray*}
which is equivalent to:
\begin{eqnarray*} & & (\exists \x_{1},...,\x_{n})(\forall \y_{1},...,\y_{m})
[(\x_{1}\in D\wedge ...\wedge\x_{n}\in D) \wedge \\ & & \wedge (\y_{1}\notin D\vee...\y_{m}\vee \io (\G)(\x_{1},...,\x_{n},\y_{1},...,\y_{m}))]\end{eqnarray*}
Since subformul{\ae}   of the form $\y \notin D$ are  negations of atomic formul{\ae} in the language of $\Q$, they are equivalent to a formula in $\Sigma _{1}^{+}$ and, by induction $\io (\G)\in \Sigma_{n-1+s}^{+}$ ($n+1$ even) or $\Sigma_{n-2+s}^{+}$ ($n+1$ odd), it follows  that the subformula 
$$(\x_{1}\in D\wedge ...\wedge\x_{n}\in D)\wedge (\y_{1}\notin D\vee...\y_{m}\vee \io (\G)(\x_{1},...,\x_{n},\y_{1},...,\y_{m}))$$
is $\Sigma _{n-1+s}^{+}$ ($n+1$ even) or $\Sigma_{n+s}^{+}$ ($n+1$ odd). Hence $\io (\Fr)$ is $\Sigma _{n+1+s}^{+}$ ($n+1$ even) or $\Sigma_{n+s}^{+}$ ($n+1$ odd). Furthermore, $t(\io(\Fr))=dm+t(\io(\G))$, from which we find by iteration that $t(\io(\Fr))=dt(\Fr)+t(\io(\times))$.\\ 
If $\Fr $ is a $\Pi_n^+$-formula, the result can be proven in a similar way -- but one has to start the induction at $n=1$ and $n=2$. 

\vv

(b)  Assume $\iota(\times) \in \Pi_s^+$. If $\Fr $ is a $\Sigma_0^+$-formula we get the same translation as in (a). The subformula
$$\bigwedge_{\bar{i}\in I}(\uu _{i_{1}},\uu _{i_{2}},\uu _{i_{3}})\in \io (+) 
 \bigwedge_{\bar{j}\in J}(\uu _{j_{1}},\uu _{j_{2}},\uu _{j_{3}})\in \io (\times)\wedge (\uu _{r},\uu _{r},\io (0))\in \io (+)) $$
is a conjunction of $\Sigma _1^+$-formul{\ae}   and $\Pi_{s}^+$-formul{\ae}, and hence is itself a $\Pi_{s}^+$-formula. The translation of $\io (\Fr )$ is then of the form:
$$(\exists \uu _1,...,\uu _l) \G$$ where $\G$ is $\Pi_s^+$, hence $\io (\Fr)$ is $\Sigma _{s+1}^+$.\\
 
If $\Fr $ is a $\Sigma_1^+$-formula, then $\Fr= (\exists x_1,...,x_n)\G(x_1,...,x_n)$ for some quantifier-free formula $\G$. The translation $\io (\Fr)$ is then given by:
 $$(\exists \x_1,...,\x_n)(\x_1\in D\wedge ...\wedge \x_n\in D \wedge \io(\G)(\x_1,...,\x_n))\; .$$ Since $\io (\G)$ is $\Sigma _{s+1}^+$ and membership of $D$ is atomic, the quantifier-free part of this formula is $\Sigma_{s+1}^+$. Hence $\io(\Fr)$ is $\Sigma_{s+1}^+$.\\
 We proceed by induction. Let $\Fr$ be a $\Sigma _{n+1}^+$-formula, i.e. 
 $$\Fr= (\exists x_1,...,x_n)(\forall  y_1,...,y_m)(\G( x_1,...,x_n,y_1,...,y_m))$$
 for some $\G \in \Sigma_{n-1}^+$. The translation $\io (\Fr)$ then becomes:
\begin{eqnarray*} & & (\exists \x_{1},...,\x_{n})(\forall \y_{1},...,\y_{m})[(\x_{1}\in D\wedge ...\wedge\x_{n}\in D)\wedge \\ & & \wedge ((\y_{m}\in D\wedge ...\wedge\y_{m}\in D)
\Rightarrow \io (\G)(\x_{1},...,\x_{n},\y_{1},...,\y_{m}))]\end{eqnarray*}
which is equivalent to:
\begin{eqnarray*} & & (\exists \x_{1},...,\x_{n})(\forall \y_{1},...,\y_{m})
[(\x_{1}\in D\wedge ...\wedge\x_{n}\in D)\wedge \\ & & \wedge (\y_{1}\notin D\vee...\y_{m}\vee \io (\G)(\x_{1},...,\x_{n},\y_{1},...,\y_{m}))]\; . \end{eqnarray*}
By induction $\io (\G)$ is a $\Sigma_{n+s}^+$-formula (if $n+1$ is even) or $\Sigma_{n+s-1}^+$-formula (if $n+1$ is odd), hence 
$$(\x_{1}\in D\wedge ...\wedge\x_{n}\in D)\wedge (\y_{1}\notin D\vee...\y_{m}\vee \io (\G)(\x_{1},...,\x_{n},\y_{1},...,\y_{m}))\; $$
is $\Sigma_{n+s}^+$ or $\Sigma_{n+s-1}^+$. From which it easily follows that $\io(\Fr)$ is $\Sigma _{n+s+2}^+$ ($n+1$ even) or $\Sigma _{n+s+1}^+$ ($n+1$ odd). \\
If $\Fr $ is a $\Pi_n^+$-formula, the result can be proved in a similar way (but starting the induction at $n=1$ and $n=2$.  \qed

\vv

\paragraph {\bf Remark.} \ If membership of $D$ is positive-existential, then one can slightly alter
the model $(D,\iota)$ to another $(D',\iota')$ in which membership of $D'$ is
quantifier-free, and hence for this altered model the theorem is true.

\vv

\paragraph {\bf Corollary.} \ \label{undecQ} {\sl If $(D,\iota)$ is a model of $\Z$ in $\Q$ that has $D$ defined by an atomic ($\Sigma_0^+$-)formula, $\iota(+)$ diophantine ($\Sigma_1^+$) and $t(\iota(\times)) \leq 1$, then the $\Sigma^+_{3}$-theory of $\Q$ is undecidable.}

\vv

\pf Davis, Matijasevich, Putnam and Robinson (cf.\ \cite{Davis:73}) have shown that the $\Sigma^+_1$-theory of $\Z$ is undecidable, but the proposition implies that any $\Sigma^+_1$-sentence 
is translated into a $\Sigma^+_{3}$-sentence over $\Q$ using $\iota$. Indeed, $\iota(\times)$ is $\Pi^+_2, \Sigma_2^+$ or $\Sigma_3^+$, and in each of these cases, a $\Sigma^+_1$-sentence translates to a $\Sigma^+_3, \Sigma_2^+$ and $\Sigma_3^+$-sentence, respectively. \qed

\vv

\vv

\sectioning{Preliminaries on elliptic divisibility sequences} \label{eds}

\vv

\paragraph {\bf Elliptic curve model of $(\Z,+)$.} \ Let $E$ denote an elliptic curve of rank one over $\Q$. Thus, as a group, $E(\Q)=\Z \oplus {\mathcal T}$ for a finite group $\mathcal T$ of cardinality $\tau$. Let $P$ be a point of infinite order on $E$. Choose a plane model $f(x,y)=0$ for $E$.

\vv

\paragraph{\bf Lemma.} \label{diodefE} \ {\sl \textup{(i)} In the above coordinates $(x,y)$, for any $r$, the set $T_{r\tau}= \langle r\tau P \rangle = \{ nr \tau P \ : \ n \in \Z \}$ is diophantine over $\Q$.

\textup{(ii)} Consider $D_r:=\{ (x,y,1) \ : \ (x,y) \in T_{r\tau} \} \cup \{ (0,1,0) \}$. Consider ${\bf 0}:=(0,1,0)$ as a symbol for the neutral element of $E$. If $+$ denotes the addition on $E$, then $(D_r,\iota)$ is a three-dimensional diophantine model of $(\Z,+)$ over $\Q$ (where $\iota(0)={\bf 0}$ and $\iota(n)=(x(n\tau rP),y(n\tau r P),1)$). Furthermore, membership of $D_r$ (``$(x,y,z) \in D_r$'') can be expressed by an atomic formula.

\textup{(iii)} The relations ``$0$'' and ``$\neq$'' in $(\Z,+)$ are diophantine over $\Q$ via $(D_r,\iota)$.}

\vv

\pf (i) Let $Q$ be a generator for the free part of $E$. Then there exists an integer $N$ such that $P=NQ$. Then $$T_{r\tau} =\{ R \in E(\Q) \ : \ (\exists S \in E(\Q))(R=Nr \tau S) \}.$$ The statement that ``$R \in E(\Q)$'' is a quantifier-free formula in $\Q$. The statement that $R=Nr \tau S$ (for fixed integers $N,r$) is too. Hence $T_r$ is diophantine over $\Q$. 

(ii) The map $\iota$ is a bijection since we have killed the torsion subgroup of $E(\Q)$ by multiplying by $\tau$. The addition formul{\ae} on $E$ can be written down in terms of coordinates on the chosen model. They will involve a choice distinction (e.g., doubling a point is different from adding two distinct points that are not opposite), but these choices are written by a formula involving inequalities and connectives, which translates into normal form only involving existential quantifiers. Hence addition is given by a diophantine formula. The statement about membership is immediate. 

(iii) $\iota(0)=(0,1,0)$ is obviously atomic. Since we are in a group, to define ``$a \neq b$'' in a diophantine way, it suffices to define ``$n \neq 0$'', and this is clearly equivalent to $\iota(n) \in T_{r\tau}$, which is diophantine. \qed

\vv

\paragraph{\bf Remark.} \ Note that if $E$ is an elliptic curve of rank one over $\Q$, there is an algorithm to compute the torsion subgroup, and if a point $P$ of infinite order is known, then one can find $N$ and $Q$ algorithmically by going through the (finite) list of points $R$ of height smaller than $P$ and checking whether $mR=P$ for the appropriate finite list of integers $m$.

\vv

\paragraph {\bf An ``odd'' divisibility sequence.} \ \label{divpoly} Let $E$ be an elliptic curve over $\Q$ of non-zero rank over $\Q$. 
Let $P$ be a point of infinite
order on $E$. We want to study arithmetical properties of the numerator and denominator of the coordinates of multiples of $P$. Choose a plane Weierstrass model for $E$:

$$ y^2 = x^3 + ax^2 + b x + c, $$
with $a,b,c$ integers. We can write
$$ nP = (x_n, y_n) = \left( \frac{a_n}{B_n^2},\frac{c_n}{B_n^3} \right),$$ with $a_n,B_n$ and $c_n,B_n$ pairs of coprime integers (with $B_n$ and $c_n$ defined up to sign). 
 
\vv

\paragraph {\bf Notation.} \ We write $(a,b)$ to denote any greatest common divisor of integers $a$ and $b$ (hence this symbol doesn't have a well-defined sign).

\vv

\paragraph {\bf Lemma.} \label{fg} {\sl \textup{(i)} If $v$ is a valuation for which $v(B_n)>0$ then for any integer $t$, $v(B_{tn}) = v(B_n) + v(t).$

\textup{(ii)} $\{ B_n \}$ is a divisibility sequence, i.e., if $m|n$, then $B_m$ divides $B_n$.

\textup{(iii)} $\{ B_n \}$ is a strong divisibility sequence, i.e., $(B_m,B_n)=B_{(m,n)}.$ }

\vv

\pf (i) For $v \neq v_2$, the claim follows from looking at the formal group law associated to $E({\bf Q}_p)$, cf.\ \cite{Cheon:98} - but some care should be taken with this reference, cf.\ the remark below. The following considerations hold regardless of the fact whether $E$ is in global minimal form or not, as long as the coefficients are integral. 

Let $v=v_p$. Let $\hat{E}$ denote the formal group of $E$. If $E_1$ denotes the kernel of reduction modulo $p$, then $E_1 \rightarrow \hat{E}(p\Z) : P=(x,y) \mapsto z(P):=-\frac{x}{y}$ is an isomorphism such that $v(z)=-\frac{1}{2} v(x)$. Theorem IV.6.4(b) from \cite{SilvermanAEC} says that if $r>1/(p-1)$, then the formal logarithm induces an isomorphism $\hat{E}(p^r\Z) \cong p^r\Z$. 

Note that for rational primes, $1/(p-1)<1$ unless $p=2$. Hence if $p \neq 2$ or $v(z(P))>1$, the isomorphism $\hat{E}(p^r\Z) \cong p^r\Z$ holds for $r=v(z(P))$ and since it is true for any larger $r$, the map preserves valuations. Hence $v(z(nP))=v(n)+v(z(P))$. This implies the claim since $z(nP)=-\frac{a_n B_n}{c_n}$ and $a_n$ and $c_n$ are coprime to $B_n$.

We are only left to consider the case $v=v_2$ and $v_2(z(P))=1$. Assume $v_2(x)<0$. The duplication formula gives 
$$ x_2 = \frac{x}{4} \cdot \frac{1-2bx^{-2}-8cx^{-3}+(b^2-4ac)x^{-4}}{1+ax^{-1}+bx^{-2}+cx^{-3}}, $$
The second factor in this product has valuation zero, and hence we get
$v(x_2)=v(x)-4$, and this implies the result for $t=2$. It follows by induction for $t=2^\ell$ for some $\ell$, and then, using the first part of the proof, for general $t=2^\ell \cdot t'$ with $t'$ odd.  

(ii) follows immediately from (i). 

For (iii), we only need to prove that $(B_m,B_n)$ divides $B_d$ for $d=(m,n)$. Choose integers $x,y$ such that $xm+yn=d$. Then part (i) implies that $v(B_{xm}) \geq v(B_m) \geq v(B_d)$ and $v(B_{yn}) \geq v(B_n) \geq v(B_d)$. Therefore $dP=xmP+ynP$ belongs to the {\sl group} of points $P=(x,y)$ with $v(x) \leq -2r$ (including the zero element of $E$) --- this is a group, since under the isomorphism $E_1 \mapsto \hat{E}(p\Z)$ it corresponds to the subgroup $\hat{E}(p^r\Z)$. Hence $v(x_d) \geq -2r$, so $v(B_d) \leq r$. \qed

\vv

\paragraph {\bf Remark.} As Marco Streng notes, the claim in \cite{Cheon:98} that $v(B_{tn}) = v(B_n) + v(t)$ also holds for the long Weierstrass form and for number fields is wrong; a counterexample is given by $P=\left(-\frac{1}{4},\frac{7}{8} \right)$ on $y^2+xy=x^3+x^2-2x$, for which $v_2(B_2)=3$ but $v_2(B_1)+v_2(2)=2$. Over an arbitrary number field, it might go wrong for a larger number of (too ramified) valuations. 

\vv

In proofs to follow, we will rely on properties of division polynomials $\phi_n, \psi_n, \omega_n$ (e.g., \cite{SilvermanAEC}  III.3.7 for standard Weierstrass form and \cite{Ayad:92} for the general case). The sequence $\{ \psi_n \}$ has been termed an {\sl elliptic divisibility sequence} by Morgan Ward (\cite{Ward:48}). This recourse to the literature is strictly speaking not necessary in this section (but we will need it in the final part of the paper), since all properties can be checked by direct, but sometimes tedious, computation using the addition formul{\ae} on $E$. Instead, we will use the following

\vv

\paragraph {\bf Substitution principle.} \ {\sl Let $f \in {\bf C}[x_1,y_1,z_1,\dots,x_r,y_r,z_r]$ be a homogeneous polynomial w.r.t.\ the weights wt$(x_i)=2i^2$, wt$(y_i)=i^2$ and wt$(z_i)=3i^2$. Suppose  $f(\phi_1,\psi_1,\omega_1,\dots,\phi_r,\psi_r,\omega_r)=0.$ Then for any point $P \in E(\Q)$ that is {\em non-singular modulo all primes}, we have
$$ f(a_1,B_1,c_1,\dots,a_r,B_r,c_r)=0,$$ if we choose the signs of $a_i,B_i,c_i$ such that they agree with those of the classical division polynomials.}

\pf The trick is dehomogeneization w.r.t.\ the denominator of $x_1$. As Mohamed Ayad has notes by direct computation in \cite{Ayad:92} (bottom of page 306), for any $n$ we can write $$x_n = \frac{b_1^{2n^2} \phi_n}{b_1^2((b_1^{n^2-1} \psi_n)^2}, y_n = \frac{b_1^{3n^2} \omega_n}{b_1^3 (b_1^{n^2-1} \psi_n)^3},$$ where numerators and denominators in these fractions are {\sl integers}; that there is no cancellation of factors of $b_1$ in this representation; and that the common divisors of $b_1^{2n^2} \phi_n$ and $b_1^{n^2-1} \psi_n$ (and $b_1^{3n^2} \omega_n$ and $b_1^{n^2-1} \psi_n$) are the primes $p$ for which $P$ is singular modulo $p$. Therefore, if $P$ is non-singular modulo all primes, we find $a_n=b_1^{2n^2}\phi_n, B_n=b_1^{n^2}\psi_n$ and $c_n=b_1^{3n^2}\omega_n$, and the result follows. \qed

Now let $E$ be an elliptic curve of rank one over $E$ with a rational two-torsion point. By translation, we can assume that $(0,0)$ is a two-torsion point on $E$. Then $E$ has a Weierstrass equation $y^2=x^3+ax^2+bx$.

\vv

\paragraph {\bf Lemma/Definition.} \label{Moh} {\sl Let $E$ be in Weierstrass form $y^2=x^3+ax^2+bx$, having $(0,0)$ as rational 2-torsion point. Let $P$ be a point of infinite order in $2E(\Q)$ (i.e., divisible by 2 in $E(\Q)$) that is non-singular modulo all primes. Then we can write 
$$ nP = (x_{n},y_{n}) = \left(\left(\frac{{A}_{n}}{B_{n}}\right)^2, \frac{A_n{C}_{n}}{B_{n}^3}\right) $$
for integers ${A}_{n},B_{n}$ and $C_n$ (defined up to sign) with $(A_n,B_n)=1$ and $(B_{n},{C}_{n})=1$. Then:

\textup{(i)} The greatest common divisor of $A_{n}$ and $C_{n}$ divides the coefficient $b$ of the Weierstrass model, and the order of $b$ at any common divisor of $A_{n}$ and $C_{n}$ is at least 2; in particular, if $b$ is squarefree, then $(A_n,C_n)=1$;

\textup{(ii)} We have $B_{2n} = {2A_nB_nC_n}$ up to sign; in particular, $A_n$ divides $B_{2n}$. }

\vv

\pf Let $P=2Q$ with $Q \in E(\Q)$. We have $\phi_2 = (\phi_1-b)^2$, so applying the substitution principle to this equation and the point $nQ$, we find that $x_n$ is a rational square and hence $A_n$ is well-defined up to sign. Substituting the point $nP$ into the equation of $E$ gives that $c_{n}^2 =A_n^2(A_{n}^4+aA_{n}^2 B_{n}^2+bB_{n}^4),$ so $A_n$ divides $c_n$ and the definition of $C_n$ makes sense. Then $C_{n}^2 =A_{n}^4+aA_{n}^2 B_{n}^2+bB_{n}^4,$ and the gcd of $C_{n}$ and $A_{n}$ divides $bB_{n}^4$, hence $b$ since $B_{n}$ and ${A}_{n}$ are coprime; if $v((C_{n},A_{n}))>0$, then we see immediately from this formula that $v(b) \geq 2$. This proves (i). 

(ii) This follows from the substitution principle via the identity of division polynomials $\psi_2 = 2 \psi_1 \omega_1$ applied to $nQ$. \qed

\vv

\paragraph{\bf Remark.} The numbers $A_n$, $B_n$ and $C_n$ (and the symbol $\sqrt{x_n}$ occasionally to be used) are only defined up to sign, but that sign will play no r\^ole in the formul{\ae}  under consideration (such as (ii) above), so we will not mention this issue anymore, except in the final section of this paper.

\vv

In subsequent considerations, we will also need to study the divisibility properties of the sequences $\{ A_\ast \}$ and $\{ C_\ast \}$. It turns out that divisibility between their $m$- and $n$-th term is only assured if $n$ is an {\sl odd} multiple of $m$. 

\vv

\paragraph{\bf Definition.} \ We call a sequence of integers $\{X_\ast\}$ an {\sl odd divisibility sequence} if $X_n$ divides $X_{nt}$ as soon as $t$ is odd. We call $\{C_\ast\}$ as defined by the previous lemma {\sl the odd divisibility sequence associated to $(E,P)$.}

\vv

That the previous definition makes sense is the contents of the following lemma:

\vv

\paragraph {\bf Lemma.} \label{odd} {\sl Assume $(E,P)$ and $(A_\ast,B_\ast,C_\ast)$ are as in Lemma \ref{Moh}, with $b$ and $a^2-4b$ squarefree. Then: 

\textup{(i)} $\{ A_\ast B_\ast \}$ is a strong divisibility sequence.

\textup{(ii)} $\{ A_\ast \}$ and $\{ C_\ast \}$ are odd divisibility sequences.

\textup{(iii)} If $t$ is odd and $v(A_n)>0$, then $v(A_{nt})=v(A_n)+v(t)$; but if $t$ is even, then $(A_n,A_{nt})=1$ for all $n$. Identical statements hold with $A_\ast$ replaced by $C_\ast$.} 

\vv

\pf  Recall that we have a morphism of 2-descent (cf.\ \cite{SilvermanAEC}, X.4.9) given by the rational map: $$\delta' \ : E \rightarrow E' \ : \ (X,Y) \mapsto (\frac{Y^2}{X^2}, \frac{Y(X^2-a^2+4b)}{X^2}) $$ with $E' \ : \ 
y^2=x^3-2ax^2+(a^2-4b)x$.

(i) Suppose $Q=nP$ maps via $\delta'$ to $Q'$. Then $x(Q')=(\frac{y(Q)}{x(Q)})^2$, so $$\sqrt{x(Q')}= \frac{A'_n}{B'_n} = \frac{C_n}{B_n A_n}$$ is a coprime representation (since $b$ is squarefree), and we find that $\{ A_\ast B_\ast \}$ is a strong divisibility sequence, as in Lemma \ref{fg}, (as it is equal to the ``$B'$''-sequence $\{B'_\ast\}$ associated to $(E',\delta'(P))$. 

(ii) Let us now prove that $\{ A_\ast \}$ is an odd divisibility sequence. Suppose $v(A_n)>0$. Then $v(B_n)=0$ by coprimeness of the representation. Now $v(A_n B_n) \geq v(A_n) >0$, and since $B'_n = 2 A_n B_n $, the formal group law on $E'$ (\ref{fg}) implies that $v(B'_{tn})=v(B'_n)+v(t)$, so we find $$ v(A_{tn})+v( B_{tn})=v(A_n)+v(t). \leqno{\indent\textrm{{\rm (\ref{odd}.1)}}}$$ If $v(B_{tn})=0$, we indeed find that $v(A_{tn}) \geq v(A_n)$. If on the other hand, $v(B_{tn})=0$, then since $A_{tn}$ and $B_{tn}$ are coprime, we find that $v(A_{tn})=0$, and hence $v(B_{tn})=v(A_n)+v(t)$.
Now $A_n$ divides $B_{2n}$ (\ref{Moh} (ii)), so we have that $(B_{tn}, B_{2n}) = B_{(tn,2n)} = B_{n(t,2)}$ is divisible by a valuation $v$ which doesn't divide
$B_n$; hence $(t,2) \neq 1$ and $t$ is even; which we have excluded. 

That $\{ C_\ast \}$ is an odd divisibility sequence is immediate, since $C_n = A'_n$ for the image sequence under $\delta'$ (with $a^2-4b$ squarefree), and we have just shown that $\{A'_\ast\}$ is an odd divisibility sequence. 

(iii) This is implicit in the proof of (ii), noting again that $A_{tn}$ and $B_{tn}$ are coprime in (\ref{odd}.1). \qed

\vv

\paragraph {\bf Remark.} Here is a quick proof that $\{A_\ast\}$ and $\{C_\ast\}$ are odd divisibility sequences: if an integer $d$ divides $A_n$ or $C_n$, then $nP \in E[2]({\bf Z}/d)$, so for odd $t$, $tnP=nP \in E[2]({\bf Z}/d)$. 

\vv

\paragraph {\bf Example.} \label{exa} \ The elliptic curve $E \ : \ y^2=x^3+12x^2+11x=x(x+1)(x+11)$ is of rank one over $\Q$, and $P=(1/4,15/8)$ is of infinite order. The torsion subgroup of $E(\Q)$ is $\Z/2 \times \Z/2$, generated by $(-1,0)$ and $(0,0)$. We computed the prime factorisations of $A_n,B_n$ and $C_n$ for $n \leq 8$:

\vv

{\small

\noindent \begin{tabular}{lcl}
$A_1$ & = &  $1$  \\
$A_2$ & = &  $5 \cdot 7$  \\
$A_3$ & = &  $19 \cdot 269 $  \\
$A_4$ & = &  $659 \cdot 1931$ \\
$A_5$ & = & $ 23042506969 $  \\
$A_6$ & = &  $\underline{5 \cdot 7} \cdot 89 \cdot 4639 \cdot 4575913 $  \\
$A_7$ & = &  $647873811 \cdot 19522768049$ \\
$A_8$ & = & $ 1321 \cdot 6637 \cdot 1356037 \cdot 6591431535431$  \\
\end{tabular}

\vv

\noindent \begin{tabular}{lcl}

 $B_1$ & = & $2$ \\
$B_2$ & = &  $\underline{2^2} \cdot 3$ \\
$B_3$ & = &  $\underline{2} \cdot 29 \cdot 41$ \\
 $B_4$ & = &  $\underline{2^3 \cdot 3} \cdot 5 \cdot 7 \cdot 37 \cdot 53 $\\
$B_5$ & = &  $\underline{2} \cdot 11 \cdot 6571 \cdot 10949$ \\
 $B_6$ & = &  $\underline{2^2 \cdot 3^2} \cdot 19 \cdot \underline{29 \cdot 41} \cdot 269 \cdot 467 \cdot 2521$ \\
 $B_7$ & = &  $\underline{2} \cdot 31 \cdot 211 \cdot 1481 \cdot 8629 \cdot 184598671 $\\
 $B_8$ & = &  $\underline{2^4 \cdot 3^1 \cdot 5 \cdot 7} \cdot 13 \cdot \underline{37 \cdot 53} \cdot 659 \cdot 1931 \cdot 160117 \cdot 5609521$ \\
\end{tabular}

\vv

\noindent \begin{tabular}{lcl}

$C_1$ & = &$ 3 \cdot 5$ \\
$C_2$ & = & $- 37 \cdot 53$ \\
$C_3$ & = &  $\underline{3^2 \cdot 5} \cdot 467 \cdot 2521$\\
$C_4$ & = &  $-13 \cdot 160117 \cdot 5609521$\\
$C_5$ & = &  $\underline{3 \cdot 5} \cdot 17 \cdot 67 \cdot 1601 \cdot 3019 \cdot  17417 \cdot 379513$ \\
$C_6$ & = &  $23 \cdot \underline{37 \cdot 53} \cdot 59 \cdot 10531 \cdot 1131223 \cdot 7186853449441$  \\
$C_7$ & = & $-\underline{3 \cdot 5} \cdot 353 \cdot 1483 \cdot 17609 \cdot 11748809 \cdot 281433601 \cdot 46333351129459$ \\
$C_8$ & = & $ 5303 \cdot 108739 \cdot 1830931 \cdot 170749043903 \cdot 92397921271034416798380481$ \\
\end{tabular}
}

\vv

\vv

Note that although $b=11$ is squarefree in the example, $a^2-4b=100$ is not. This means something might go wrong with the valuation formula for $C_\ast$ upon multiplication by 2 or 5, and indeed, $v_5(C_5) \neq v_5(C_1)+v_5(5)$.

The examples illustrate all the (non-)divisibility-properties mentioned before, but also some other apparent features that will be discussed later on: whereas the indices have one prime factor on average, the numbers themselves have three primitive factors on average. It is expected that for any given $k>0$, all terms in the sequence from a certain moment on will have at least $k$ primitive factors. In the above tables, we have underlined the ``non-primitive'' part, i.e., the prime factors that occur earlier on the list. 

\vv

\indent {\bf Observation.} \label{rubin} \ {\sl All divisors of $A_n$ and $B_n$ for odd $n$ are $= \pm 1$ mod 5. }

\vv

{\indent\textit{Proof} (for $B_\ast$, as shown to us by Karl Rubin).\ } Suppose $l$ is a prime with $l|B_n$, i.e., $nP=0$ mod $l$. Since $n$ is odd, $P=2Q$ mod $l$ for $Q=(n+1)/2 \cdot P$. Then $x=x(Q)$ satisfies the equation $(x^2-8x+11)(x^2+7x+11)=0$ mod $l$. Since both factors of this equation have discriminant $5$ up to squares, there is a solution mod $l$ precisely if $5$ is a square modulo $l$. \hfill $\Box$

\vv

On the other hand, all $C_n$ seem to have primitive prime divisors of odd order $=\pm 2$ mod 5, i.e., inert in $\Q(\sqrt{5})$, but we have no general proof of that.

\vv

\vv

\sectioning{Elliptic divisibility sequences and models of $(\Z,+,|)$} \label{edsmod}

\vv

\paragraph {\bf Primitivity condition.} \ \label{primcond} Let $\{X_\ast\}$ be an (odd) divisibility sequence. Let $R$ denote a set of valuations. We say $\{X_\ast\}$ {\sl is $R$-primitive} if every term $X_n$ has a primitive divisor from $R$, that is: 
$$ (\forall n)(\exists v \in R)[v(X_n)>0 \mbox{ and } (\forall i<n) (v(X_i)=0)]. $$
We say $\{X_\ast\}$ {\sl is $R$-odd-primitive} if every term $X_n$ has a primitive {\sl odd order} divisor from $R$, that is: 
$$ (\forall n)(\exists v \in R)[v(X_n) \mbox{ is odd and } (\forall i<n) (v(X_i)=0)]. $$ We sometimes say $v$ is $R$-(odd-)primitive for $X_n$ if these formul{\ae}  holds for $v$ and $X_n$. 

\vv

\paragraph {\bf Lemma.} \label{Cdiv} \ {\sl Suppose that $E$ and $P$ are as in lemma \ref{Moh}. Assume that $\{ C_\ast \}$ is $R$-(odd-)primitive for some $R$. If $v \in R$ is (odd-)primitive for $C_m$ and $v(C_n)>0$ for some $n$, then $m|n$ and $n/m$ is odd.}

\vv

\pf  It suffices to prove this for the $A$-sequence, since the descent morphism $\delta'$ transfers $\{C_\ast\}$ into $\{A_\ast\}$ (proof of Lemma \ref{odd}). Now since $\{ A_\ast B_\ast \}$ is a strong divisibility sequence (\ref{odd}), we have 
$$ (A_m B_m, A_n B_n) = A_{(m,n)} B_{(m,n)}. $$
Since we assume $v(A_m)>0$ and $v(A_n)>0$, we have $v(B_m)=v(B_n)=0$ by coprimeness assumptions; and $v( (A_{(m,n)} B_{(m,n)}))>0$ by the above formula. 
Suppose first that $v(A_{(m,n)})>0$. Since $(m,n) \leq m$, the $R$-primitivity of $v$ for $A_m$ implies that $(m,n)=m$. This means that $m|n$. By \ref{odd} (iii), we find that $n/m$ is odd. 

On the other hand, if $v(B_{(m,n)})>0$, since $\{ B_\ast \}$ is a divisibility sequence and $(m,n)|m$, we have $v(B_m)>0$, contrary to the assumption. \qed

\vv

\paragraph {\bf Divisibility predicate.} \ Let $R$ denote a set of valuations. Denote by $\mathcal{D}_R(x,y)$ the property
$$\mathcal{D}_R(x,y) \ : \ \forall v \in R \ : \ v(x) \mbox{ odd } \Rightarrow v(x) < v(y^2). $$

\vv

\paragraph {\bf Remark.} \ This predicate says that odd order ``zeros'' of $x$ are zeros of at least half that order of $y$, and that odd order ``poles'' of $x$ are at most poles of $y$ of half that order. Note that it seems at this point maybe more natural to have a definition in which the condition $v(x)<v(y^2)$ is replaced by $v(x)<v(y)$, but for future applications, we will need it as it stands.

\vv

\paragraph {\bf Theorem.} \label{defdiv} {\sl Let $E$ be an elliptic curve over $\Q$ and $P$ a point of infinite order on $2E(\Q)$ of sufficiently large height. Assume $E$ has Weierstrass form $y^2=x^3+ax^2+bx$ (in particular, a rational 2-torsion point) with $b$ and $a^2-4b$ squarefree. Assume the odd divisibility sequence $\{C_\ast\}$ associated to $P$ on $E$ is $R$-odd-primitive. Then for any integers $m,n \in {\bf Z}$,
\begin{eqnarray*} 
m | n  &\iff& \D_R(y_m \sqrt{x_m}, y_n \sqrt{x_n}) \vee \D_R(y_m \sqrt{x_m}, y_{m+n} \sqrt{x_{m+n}})
\end{eqnarray*}}
\pf Replacing $P$ by a suitable multiple, we can assume $P$ is non-singular modulo all primes. Indeed, for any prime $p$, consider the group $E(\Q_p)$ and the subgroup $E_0(\Q_p)$ of points that reduce to non-singular points modulo $p$. Then $E(\Q_p)/E_0(\Q_p)$ is finite and non-zero for only finitely many $p$ (actually, by a theorem of Kodaira and N\'eron, of order bounded uniformly in $p$ by 4 times the least common multiple of the exponents in the minimal discriminant of $E$, cf. \cite{SilvermanAEC} VII.6.1).  Note that the $R$-odd-primitivity condition is unaffected by this replacement of $P$ by a multiple. By \ref{odd}, $\{C_\ast\}$ is an odd divisibility sequence. 

It follows from the definition of $C_n$ that 
$$ C_N = \pm y_N \sqrt{x_N} \left(\frac{B_N^2}{A_N}\right)^2. \leqno{\indent\textrm{{\rm (\ref{defdiv}.1)}}}$$
We claim that our assumption that $b$ is squarefree implies the following:
 
\vv

{\bf \indent\textbf{{\ref{defdiv}.2}} Claim.} \ {\sl  If $v(C_N) \neq 0$, then $v(C_N)=v(y_N \sqrt{x_N})$.}

\vv

{\sl Proof of the claim.}\ By the above formula we should prove $v(B_N^2/A_N)=0$. Now $B_N$ and $C_N$ are coprime by definition, and by \ref{Moh} (i) and since $b$ is squarefree, we find that $A_N$ and $C_N$ are also coprime. 

\medskip

{\sl Proof of $\Rightarrow$.} \ Assume $m|n$. Then either $n/m$ or $(n+m)/m$ is odd. Then lemma \ref{odd} implies that $C_m | C_n$ or $C_m | C_{m+n}$. We will agree from now on to write $n$ but mean either $n$ or $m+n$, and assume that $n/m$ is odd. 

Pick a valuation $v$ in $R$ and suppose that $v(y_m \sqrt{x_m})$ is odd.  
From formula (\ref{defdiv}.1), we see that $v(C_m)$ has to be odd. Since $n/m$ is odd, lemma \ref{odd} implies that $v(C_n) \geq v(C_m)>0$. By (\ref{defdiv}.2), we find $v(y_n \sqrt{x_n}) \geq v(y_m \sqrt{x_m})>0$ and this implies $\D_R(y_m \sqrt{x_m}, y_n \sqrt{x_n})$. 
\medskip

{\sl Proof of $\Leftarrow$.} \ Choose a valuation $v$ that belongs to an odd order primitive divisor of $C_m$ from $R$. Then claim (\ref{defdiv}.2) implies that $v(y_m\sqrt{x_m})$ is positive and odd. The assumption means that $2v(y_n \sqrt{x_n}) > v(y_m \sqrt{x_m}) > 0$ (or similarly with $n$ replaced by $m+n$). Formula (\ref{defdiv}.1) implies that one of the following two cases has to occur: $v(C_n)>0$ or $v(A_n)>0$. In the first case, since $v$ is primitive for $C_m$, we find that $m|n$ from Lemma \ref{Cdiv}. In the second case, note that $A_n$ divides $B_{2n}$ (\ref{Moh}(ii)), so $B_{2n}$ and $C_m$ have a common divisor $v$. We will prove that $v$ is a primitive divisor of $B_{2m}$. Since $v(B_{2n})>0$, we will find from this that $m|n$. 

By \ref{Moh} (ii), we have an identity
$$ C_m = \frac{B_{2m}}{2A_m B_m}.  \leqno{\indent\textrm{{\rm (\ref{defdiv}.3)}}} $$
Should $v(B_m)>0$, then $v(B_{2m})=v(B_m)+v(2)$, and hence (as $A_m$ and $B_m$ are coprime) $v(C_m)=0$, a contradiction. Hence $$  v(B_m)=0 \mbox{ and } v(B_{2m})>v(2A_m). \leqno{\indent\textrm{{\rm (\ref{defdiv}.4)}}} $$
Since $v$ is primitive for $C_m$, we have $v(C_i)=0$ if $i|m, \ i<m$. Recall $v(B_m)=0$, and since $\{B_\ast\}$ is a divisibility sequence, $v(B_i)=0$ for all $i$ as before. Hence
for such $i$, we find from (\ref{defdiv}.3):
$$  i|m, \ i<m  \Rightarrow v(B_{2i})=v(2A_i).  \leqno{\indent\textrm{{\rm (\ref{defdiv}.5)}}}$$

Suppose that $m/i$ is odd. Since $\{A_\ast\}$ forms an odd divisibility sequence (by \ref{odd}), should $v(A_i)>0$, then $v(A_m)>0$. But since we assume $b$ squarefree, $A_m$ and $C_m$ are coprime, so we cannot have $v(A_m)>0$ and $v(C_m)>0$. Hence $v(A_i)=0$ for all $v$ and so by (\ref{defdiv}.5), $v(B_{2i})=0$ for all $i|m$, unless $v=v_2$. For $v=v_2$, we find instead that $v(B_{2i})=v(2A_i)=1$, and hence $v(B_{2m})=v(B_{\frac{m}{i} \cdot 2i}) = v(B_{2i})=1$ by the formal group law, since $m/i$ is odd. But since $v(C_m)>0$, we have $v(B_{2m})>1$ from (\ref{defdiv}.3). This is a contradiction.

For the general case ($m/i$ not odd), write $m/i = 2^l \cdot k$ with $k$ odd. We conclude from the previous reasoning that $v(B_{2^{l+1} i})=0$, but since $\{B_\ast\}$ forms a divisibility sequence, we find from this that also $v(B_{2i})=0$. 

\vv

Recall that $v(B_m)=0$ and hence $v(B_i)=0$ for all $i|m$. We conclude from this and $v(B_{2i})=0$ that $v$ is also primitive for $B_{2m}$. But remember we had 
$v(B_{2n})>0$. Hence $2m|2n$. \qed

\vv

\paragraph {\bf Remarks.} \ (i) It might be possible to remove the assumption that $b$ or $a^2-4b$ be squarefree, but then (\ref{defdiv}.2) changes.

\vv

(ii) We work with $\{ C_\ast \}$ because $C_n$ is an algebraic function of the coordinates of $nP$ ``up to squares'' (\ref{defdiv}.1). It has been suggested in the past that $m|n$ is equivalent to $\D_R(\sqrt{x_m},\sqrt{x_n})$, but this is wrong in two ways: if $n/m$ is even, then ``zeros'' of $x_m$ are not zeros of $x_n$; and in general, ``poles'' of $x_m$ are ``poles'' of $x_n$ of {\sl larger} order (in particular, if $n/m$ is divisible by that pole). In case of an isotrivial elliptic curve over a rational function field, this problem doesn't occur (\cite{Pheidas}, esp.\ 2.2), since the order stays equal.  

\vv

Not to obscure the proof too much, we have not included the following stronger statement in the original statement of the theorem:

\vv

\paragraph {\bf Proposition/Definition.} \ \label{weakprim} {\sl The conclusion of the above theorem \ref{defdiv} still holds if one replaces the $R$-odd-primitivity condition of $\{ C_\ast \}$ by the following weaker condition: \begin{quote} \textup{(``{\bf weak} $R$-odd-primitivity condition'')} \  All terms $C_{2^a p^b}$ for $a,b$ positive integers and $p$ any odd prime have a primitive odd order divisor from $R$.  \end{quote}}
\pf \ The only part of the proof that changes is the proof of $\Leftarrow$. Set $a=v_2(m)$. If $p$ is an odd prime such that $v_p(m)=b>0$, choose a primitive odd divisor $v$ for $C_{2^a p^b}$ from $R$ based on the assumption. Since $m/2^a p^b$ is odd, lemma \ref{odd} implies that $v(C_m)$ is odd, so the assumption of the theorem assures us that $v(C_n)>0$. We can then proceed as before with $m$ replaced by $2^a p^b$ to conclude $2^a p^b |n$. Since this holds for any odd $p$, we find $m|n$. \qed

\vv

We now suggest the following conjecture about the sequences $\{C_\ast\}$: 

\vv

\paragraph {\bf Conjecture.} \label{GC} \ {\sl The following exist:  

\textup{(a)} an elliptic curve over $\Q$, such that $E$ has Weierstrass form $y^2=x^3+ax^2+bx$ (in particular, a rational 2-torsion point) with $b$ and $a^2-4b$ squarefree; 

\textup{(b)} a point $P$ of infinite order on $E$ with associated odd divisibility sequence $\{C_\ast\}$;

\textup{(c)} a set $R$ of prime numbers such that ${\mathcal D}_R$ is diophantine over $\Q$; and such that $\{ C_\ast \}$ is (weakly) $R$-odd-primitive.}

\vv

\paragraph {\bf Theorem.} \label{thmGC} {\sl Assume Conjecture \ref{GC}. Then $(\Z,+,|)$ has a three-dimensional diophantine model in $\Q$.}

\vv

\pf \ Immediate from \ref{diodefE}, \ref{defdiv} and \ref{weakprim}, observing that $a=\pm y_n \sqrt{x_n}$ for  $(x_n,y_n) \in T_r$ is diophantine over $\Q$. \qed

\vv

\paragraph {\bf Proposition (($C$)-elliptic Zsigmondy's theorem).} \ {\sl Let $E$ be an elliptic curve over $\Q$ and $P$ a point of infinite order in $E(\Q)$ of sufficiently large height. Let $R$ denote the set of all finite valuations of $\Q$. Then 

\textup{(i)} $\{B_\ast\}$  is $R$-primitive. 

\textup{(ii)} If $(0,0) \in E[2]$, then $\{ C_\ast \}$ is $R$-primitive.

\textup{(iii)} If $E$ has $j$-invariant $j=0$ or $j=1728$, then the $ABC$-conjecture implies that $\{B_\ast\}$ and $\{ C_\ast \}$ (for $(0,0) \in E[2]$) are $R$-odd-primitive.}

\vv

\pf The crucial statement is Siegel's theorem on integral points on an elliptic curve (cf.\ \cite{SilvermanAEC}, IX 3.3), which implies that $A_n$ and $B_n$ are both of order of magnitude the height of $nP$. 

For $B_\ast$, (i) is the usual elliptic Zsigmondy's theorem, first proven by Silverman in \cite{SilvWief}, Lemma 9. The same proof works for the sequence $\{ A_* \}$; we include a variation of the proof for completeness.

Suppose that $A_n$ doesn't have a primitive divisor. We will show that $n$ is absolutely bounded, so changing $P$ to some multiple, we get the result. We claim that there exists a set $W$ of distinct divisors $d$ of $n$ with all $d>1$ such that 
$$ A_n \ | \ n \prod_{ d \in W} A_{\frac{n}{d}}. $$
We can then finish the proof as follows: We get
$$ \log |A_n| \leq \log n + \sum_{d \in W} \log |A_{\frac{n}{d}}|.$$
Let $m$ denote the canonical height of $P$. Classical height estimates give $\log |A_{\frac{n}{d}}| \leq \left( \frac{n}{d} \right)^2 m + O(1)$. They combine with Siegel's theorem (``$|A_n|$ and $|B_n|$ are of the same size'') to give for any $\varepsilon>0$, $(1-\varepsilon)n^2 m \leq \log |A_n|$. Since $${\displaystyle \sum\limits_{d \in W} \frac{1}{d^2} < \zeta(2)-1} \mbox{ (recall: $d>1$ for $d \in W$),} $$ we find after insertion of these estimates into the above formula:
$$ (2-\varepsilon -\zeta(2)) m n^2  \leq \log(n) + O(1), $$
and this bounds $n$ absolutely.

For the proof of the claim: by assumption, any prime $p$ dividing $A_n$ divides $A_m$ for some $m<n$. Then also $p|A_{(m,n)}$ (as in the proof of \ref{Cdiv}), so we can assume $m|n$ and $n/m$ odd. Hence $v_p(A_n)=v_p(A_m)+v_p(n/m)$ (\ref{odd} (iii)). Run through all $p$ in this way, and pick such an $m$ for each $p$. If $v_p(A_n)=v_p(A_m)$, then let $d:=n/m \in W$. Then $v_p(A_n)=v_p(A_{\frac{n}{d}})$. If, on the other hand, $v_p(A_n)>v_p(A_m)$, then we must have $p|\frac{n}{m}$. In this case, set $d:=p \in W$. Then $v_p(A_n) = v_p(A_{\frac{n}{d}})+1$. Indeed, we only have to prove that $v_p(A_{\frac{n}{d}})>0$ since we get the implication by the formal group law formula as $p$ is odd. Now since $v_p(A_m)>0$ and $n/mp$ is an odd integer, the same formula implies that $v_p(A_{\frac{n}{p}})=v_p(A_m)+v_{p(\frac{n}{mp}})>0$, and we are done. 

To finish the proof of the proposition, the statement is true for $C_\ast$, since it is the $A_\ast$-sequence of the isogenous curve $E'$ (as observed before). We note that (iii) for $\{ B_\ast \}$ is Lemma 13 in \cite{SilvWief}, and a similar argument works for $\{ C_\ast \}$. \qed

\vv

\paragraph{\bf Remarks.} \ (i) We don't know whether ${\mathcal D}_R$ for $R$ the full set of valuations is diophantine over $\Q$. This would be quite a strong statement. For example, if we write a rational number $x$ as $x=x_0 \cdot x_1^2$ with for any $v \in R$ such that $v(x_0) \neq 0$, one has $v(x_0)$ odd and $v(x_1)=0$, then ${\mathcal D}_R(x^{-1},1)$ expresses that $x_0$ is an integer.

\vv

(ii) Using elliptic Zsigmondy and a proof similar to (but easier than) that of theorem \ref{defdiv}, one can prove that $m|n \iff B_m | B_n \iff \mbox{rad}(B_m)|B_n$. However, we don't know that the formula $\Fr(x,y) : (\forall v)(v(x)<0 \Rightarrow v(y) < 0)$ is equivalent to a diophantine formula $\mathcal{D}'(x,y)$ in $\Q$. If so, then $\mathcal{D}'(x_m,x_n)$ would be a diophantine definition of $m|n$ in $\Q$. But then again, ``$\mathcal{D}'(x,1)$'' would be a diophantine definition of $\Z$ in $\Q$. 

\vv

\paragraph {\bf A diophantine predicate.} \ We will now investigate in how far one can construct sets $R$ for which ${\mathcal D}_R$ is diophantine over $\Q$. In \cite{Pheidas}, Pheidas has produced a diophantine definition over $\Q$ that says of two rational numbers $x$ and $y$ that for any prime $p = 3$ mod 4, we have $v_p(x) > v_p(y^2)$ (and  some extra conditions). This was consequently extended to all primes inert in a given quadratic extension of $\Q$ by Van Geel and Zahidi at Oberwolfach (\cite{ZahidiVanGeel}), but still involving extra conditions. Finally, Demeyer and Van Geel have proven the following (for an arbitrary extension of global fields, but we only state it for $\Q$):

\vv

\paragraph {\bf Proposition.} (\cite{Demeyer:04}) {\sl For a non-square $d$, let $R_d$ denote the set of valuations of $\Q$ that are inert in $\Q(\sqrt{d})$. Then there is a (diophantine) $\Sigma_1^+$-formula equivalent to $\mathcal{D}_{R_d}(x,y)$, i.e., $t(\mathcal{D}_{R_d})=0$.}

\vv

\paragraph {\bf Remarks.} The proof involves the theory of quadratic forms and is very close in spirit to the proof of Pheidas, which in its turn is an attempt to analyse Julia Robinson's definition $\mathcal R$ from the following perspective: $\mathcal R$ is essential a conjunction over all valuations $v$ of a predicate that says that a rational number $x$ is $v$-integral. The latter is expressed by the isotropy of a quaternary quadratic form that depends on $x$ and $v$. Pheidas' analysis says that one can discard this conjunction over an infinite set of primes (but not all). It would be interesting to see whether $\mathcal{D}_{R}$ is diophantine for other sets of primes $R$ that are inert in not necessarily quadratic extensions of $\Q$. Note that one can define $v(x) \geq 0$ for all $v$ not completely split in a cyclic extension of $\Q$ of degree $q$ (with finitely many exceptions on $v$), but for $x \in {\bf Z}[T^{-1}]$ where $T$ is the complement of finitely many primes, instead of $x \in {\bf Q}$ (Shlapentokh \cite{Shlap}, 4.4.6). 

\vv

\paragraph {\bf Corollary.} \label{dir} {\sl For any finite set $D$ of fundamental discriminants, set $$R_D := \bigcup\limits_{d \in D} R_{d}.$$ Then $\mathcal{D}_{R_D}$ is expressible by a $\Sigma_1^+$-formula. In particular, there are sets of primes $R$ of arbitrary high Dirichlet density $<1$ for which $\mathcal{D}_R$ is diophantine}.

\vv

\pf \ The first claim is automatic, since a finite disjunction of $\Sigma_1^+$-formul{\ae} is equivalent to a $\Sigma_1^+$-formula. For the second statement, choose all $\Q(\sqrt{d})$ for $d \in D$ linearly disjoint, then
$R_D$ is the complement of the set of primes that split completely in the 
compositum $L$ of all $\Q(\sqrt{d})$ for $d \in D$, and this complement has Dirichlet density $1/|L|=1/2^{|D|}$ (by Chebotarev's or weaker density theorems), which can be made arbitrary small $\neq 0$ by increasing $|D|$. \qed

\vv

Based on this information, we change our conjecture to the following, the plausibility of which will be discussed in another section, and that will be used here as input for our main theorem. 

\vv

\paragraph {\bf Conjecture.} \label{SC} \ {\sl The following exist:  

(a) an elliptic curve over $\Q$, such that $E$ has Weierstrass form $y^2=x^3+ax^2+bx$ (in particular, a rational 2-torsion point) with $b$ squarefree; 

(b) a point $P$ of infinite order on $E$ with associated odd divisibility sequence $\{C_\ast\}$;

(c) a finite set $D$ of quadratic discriminants such that $\{ C_\ast \}$ is (weakly) $R_D$-odd-primitive.}

\vv

As before, we get:

\vv

\paragraph {\bf Theorem.} \label{thmSC} {\sl Assume Conjecture \ref{SC}. Then $(\Z,+,|)$ has a diophantine model in $\Q$. \qed}

\vv

\vv

\sectioning{Defining multiplication in $(\Z,+,|)$} \label{x}

\vv

\paragraph {\bf Lemma.}\ \label{lemdiv} {\sl There exists a $\Sigma_3^+$-formula $\Fr$ in $(\Z,+,|,\neq)$ such that for integers $m,n,k$, we have $k=m \cdot n \iff \Fr(m,n,k)$.} 

\vv

\pf\ The first part of the proof is very similar to that of Lipshitz for ${\bf N}$ in \cite{Lip3}. To define multiplication by a $\Sigma_3^+$-formula, it suffices to define squaring by a $\Pi_2^+$-formula, since 
$x=mn \iff 2x = (m+n)^2-m^2-n^2$ (translating this as $(\exists u,v,w,s)(x+x+u+v=w \wedge u=m^2 \wedge v=n^2 \wedge w=s^2 \wedge s=m+n)$). 

We first claim that $y=x^2$ if and only if $(\forall t) (\phi (x,y,t))$, where $\phi$ is the formula

\begin{eqnarray*} \phi & : & x|y\wedge x+1|y+x\wedge x-1|y-x\wedge \\ & & \left(( x|t\wedge x+1|t+x\wedge 
x-1|t-x)\Rightarrow (y+x|t+x\wedge y-x|t-x)\right )\end{eqnarray*}
Indeed, the first three divisibilities imply $y=u x$ 
with $u \in {\bf Z}$ and $|x+1|\leq|u+1|$ and $|x-1|\leq|u-1|$. 
Taking $t=x^{2}$, the divisibilities following the implication sign imply 
$|x+1|\geq|u+1|$ and $|x-1|\geq|u-1|$. Hence 
$x+1=\pm (u +1)$ and $x-1=\pm (u-1)$. If in either of the two equalities, the equality 
holds with a positive sign we get that $u=x$ and hence $y=x^2$. The case
$x-1=-u+1$ and $x+1=-u-1$ leads to a contradiction. 
The other direction is easy.

Rewriting the formula $\phi$ as an atomic formula using the recipe from Lemma \ref{prenex}, we see that the replacement of the implication in $\phi$ by disjunction introduces (non-positive) expressions of the form ``$a$ does not divide $b$''. We will now show how to replace this by a positive existential statement in $(\Z,+,|,\neq)$. 

Observe that $g$ is a greatest common divisor of $a$ and $b$ in $\Z$ (notation: $(g)=(a,b)$) if and only if $$g|a \wedge g|b \wedge (\exists x,y)(a|x \wedge b|x \wedge g=x+y). \leqno{\indent\textrm{{\rm (\ref{lemdiv}.1)}}}$$
Indeed, the first two divisibilities imply an inclusion of ideals $(a,b) \subseteq (g)$, and the existential statement implies that $g \in (a,b)$.

Now $a$ doesn't divide $b$ if and only if $(a) \neq (a;b)$, and this can be rewritten as 
$$ (\exists g,g')((g)=(a,b) \wedge g+g'=0 \wedge a \neq g \wedge a \neq g'), $$
which is a positive existential formula in $(\Z,+,|,\neq)$, after substitution of (\ref{lemdiv}.1). \qed

\vv

\paragraph {\bf Theorem.} \label{42} {\sl Assume Conjecture \ref{GC} or \ref{SC}. Then $\Z$ has a model $D$ in $\Q$ with complexity $t(D)\leq 1$, $c(D)\leq 1$; and the $\Sigma_3^+$-theory ($=\Sigma_3$-theory) of $\Q$ is undecidable.}

\vv

\pf \ Picking an elliptic curve as in one of the conjectures, we find a three-dimensional diophantine model of $(\Z,+,|)$ in $(\Q,+,\times)$ as in Theorem  \ref{thmGC} or \ref{thmSC}. Now observe that $0$ is also definable in the model by an atomic formula, and that $n \neq 0$ is also definable in the model by an existential formula, cf.\ Lemma \ref{diodefE}. 
Hence each of the conjectures actually imply that $(\Z,+,|,0,\neq)$ is definable in $\Q$.

Now $\times$ is defined by a $\Sigma^+_3$-formula
 in $(\Z,+,|,\neq)$ with only one universal quantifier; in particular, $t(\iota(\times)) \leq 1$ and $c(\iota(\times)) \leq 1$ for the induced model $D$ of $\Z$ in $\Q$.  By \ref{undecQ}, we conclude that $t(D)\leq 1$, $c(D)\leq 1$ and that the $\Sigma_3^+$-theory of $\Q$ is undecidable. \qed 

\vv

\paragraph {\bf Remark.} \ The trick of replacing non-divisibilities by existential sentences in the lemma (communicated to us by Pheidas) is crucial. The negation of a diophantine formula expressing divisibility (as it comes out of our conjecture) is a universal formula that leads to a model of the same complexity as Julia Robinson's. 



\vv

\vv

\sectioning{Conditional undecidability of the $\Pi_2^+$-theory of $\Q$} \label{pi2}
\vv

Theorem \ref{42} is our main result about the complexity of a model of $\Z$ in $\Q$. Although the model given there is $\Sigma_3^+$, we can slightly alter the method to (conditionally) prove the existence of undecidable formul{\ae} in $\Q$ of complexity $\Pi_2^+$.

\vv

\paragraph {\bf Lemma}. \label{pi-undec} {\sl
The $\Pi_2^+$-theory of $(\Z,+,|)$ is undecidable.}

\vv

\pf The proof is entirely analogous to that for ${\bf N}$ by Lipshitz in \cite{Lip3}. By Lemma \ref{lemdiv}, we know that squaring is definable in $({\mathbf
Z},+,|)$ by a
$\Pi_2^+$-formula. The proof of \ref{lemdiv} actually shows that if one
allows negated
divisibilities, the defining formula can be taken to be a 
$\Pi_1$-formula.
It is easily seen that the $\Sigma_1^+$-theory of the structure $({\mathbf
Z},+,x\rightarrow x^2)$ is
undecidable (since multiplication is existentially definable). We obtain that the $\Sigma_2$-theory of $({\mathbf Z},+,|)$ is
undecidable. Since
the negation of a $\Sigma_2$-formula is a $\Pi_2$-formula, we obtain that the
$\Pi_2$-theory of $({\mathbf Z},+,|)$ is undecidable. This means that sentences of
the form:
$$\forall {\mathbf x}\exists {\mathbf y} \phi({\mathbf x},{\mathbf y})$$
(with $\phi $ quantifier free) are undecidable. However $\phi $ may still
contain negated
divisibilities and inequations. These can be eliminated as in the proof of
\ref{lemdiv} at
the expense of introducing extra existential quantifiers. Thus, any
$\Pi_2$-sentence is
equivalent to a $\Pi_2^+$-sentence and hence the $\Pi_2^+$-theory of
$({\mathbf Z},+,|)$ is
undecidable.\qed

\vv

We now need the following extension of Lemma \ref{compl-change}:

\vv

\paragraph {\bf Proposition.} \ \label{compl-change2} {\sl Let $(D,\io)$
be a diophantine model of
$(\Z,+,|)$ in $(\Q,+,\times)$, such that membership of $D$ is quantifier-free. In the
following table, the
second column lists the positive hierarchical status of the formula
$\io(\Fr)$ as a
function of the status of $\Fr$:

\renewcommand{\arraystretch}{1.2}

\begin{center}
\begin{tabular}{l|l}
 $\Fr $& $\io(\Fr)$  \\
\hline
$\Sigma_{2n}^{+}$ & $\Sigma_{2n+1}^{+}$  \\
$\Sigma_{2n+1}^{+}$ & $\Sigma_{2n+1}^{+}$  \\
 $\Pi_{2n}^{+}\; (n>0)$ & $\Pi_{2n}^{+}$  \\
 $\Pi_{2n+1}^{+}$ & $\Pi_{2n+2}^{+}$  \\
\end{tabular}
\end{center}
(note: inclusion of a formula in a class of the hierarchy means that the formula is equivalent to a formula in that class). 
} 

\vv

\pf The proof is completely analogous to the proof \ref{compl-change}. \qed

\vv

\paragraph {\bf Theorem.} \ \label{undecQ2} {\sl Assume Conjecture \ref{GC} or \ref{SC}. Then the $\Pi_2^+$-theory of
$\Q$ is
undecidable. Furthermore the subset of sentences of $\Pi_2^+$ with
$t$-complexity $1$ is
already undecidable.}

\vv

\pf  Follows immediately from  \ref{thmGC}, \ref{thmSC}, \ref{pi-undec} and \ref{compl-change2}. \qed

\vv

\vv

\sectioning{Discussion of the conjecture} \label{conj}

\vv

\paragraph {\bf Different versions of the conjecture.} \ Our conjecture is merely about the existence of {\sl one} elliptic curve, but one can of course also investigate whether the conjecture might be true for any elliptic curve $E$ with a point of infinite order on it. The conjecture then becomes a kind of elliptic Zsigmondy conjecture with odd order and inertial conditions. It then seems natural to also look at the conjecture for $\{ B_\ast \}$ instead of $\{ C_\ast \}$, although we don't know of a direct application to logic. We now first list these variants of the conjecture in a more precise way:

\vv

\paragraph {\bf (Odd-)inertial $C$-elliptic Zsigmondy's conjecture.} \ \label{ICZ} {\sl For every elliptic curve $E$ in Weierstrass form such that $(0,0) \in E[2]$ and every rational point $P$ of infinite order and sufficiently large height, the associated odd divisibility sequence $\{ C_\ast \}$ is (weakly) $R_D$-(odd-)primitive for some $D$.} 

\vv

\paragraph {\bf (Odd-)inertial elliptic Zsigmondy's conjecture.} \ \label{IZ} {\sl For every elliptic curve $E$ in generalised Weierstrass form and every rational point $P$ of infinite order and sufficiently large height, the associated elliptic divisibility sequence $\{ B_\ast \}$ is (weakly) $R_D$-(odd)-primitive for some $D$.} 

\vv

It is hard to falsify these conjectures, because if one finds a multiple of a given point $P$ for which the divisibility sequences under consideration has {\sl no} primitive odd order divisor from a given $R_D$, one simply takes a multiple of $P$ or enlarges the set of discriminants. But if the height of $P$ becomes too large, one can no longer factor $B_n$ or $C_n$ in reasonable time with existing algorithms, and if the height of $P$ is too small, then $B_n$ or $C_n$ could be non-typical (e.g., prime) for small $n$ (similar problems occur in \cite{EverestWard}). We will therefore refrain from presenting extensive numerical computation, but rather present some heuristics and remarks below, and a density version of the conjecture in the next section.

\vv

\paragraph{\bf Heuristic arguments.} \ \label{heur} 

\vv

We start from the following observation: 

\vv

(\ref{heur}.1) (Landau-Serre \cite{Serre} 2.8) \ {\sl Let $M$ be a multiplicative set of positive non-zero integers (i.e., $xy \in M \iff x \in M \vee y \in M$), and assume that the set of prime numbers in $M$ is frobenian with density $\delta>0$ (i.e., every sufficiently large prime $p$ belongs to $M$ exactly if its Frobenius morphism belongs to a fixed subset $H$ of the Galois group $G$ of some fixed number field with $H$ stable under conjugation by $G$ and $\delta=|H|/|G|$). Then the probability that a given number $x$ belongs to the complement of $M$ admits
an asymptotic expansion
$$ \log(x)^{-\delta} (\sum_{i=0}^N c_i \log(x)^{-i} +O(\log(x)^{-(N+1)})) $$
with $c_0>0$, for any positive integer $N$.}

\vv

We can now ``prove'' heuristically:

\vv

(\ref{heur}.2) \ {\sl Let $E$ be an elliptic curve over $\Q$. Then set $A=E(\Q)-E(\Z[\frac{1}{R_D}])$ of points whose denominators are only divisible by primes outside $R_D$ is heuristically finite if $|D|$ is large enough.}

\vv

Let $M$ denote the set of integers having at least one factor from $R_D$. Then $M$ is multiplicative, and a prime $p$ belongs to $M$ exactly if $p$ is not completely split in the compositum $L=\Q(\sqrt{d_1},\dots\sqrt{d_N})$, where $D=\{d_1,\dots,d_N\}$. This is the same as saying that the Frobenius element of $p$ belongs to $H=\mbox{Gal}(L/\Q)-\{1\}$. Note that $H$ is stable under conjugation, and that $\delta:=|G|/|H|=1-1/2^N>0$.

We approximate the probability that a number is outside $M$ by the first order term in (\ref{heur}.1) --- in the considerations below, any finite order truncation actually gives the same result. We consider the set $$A_x=\{ P \in E(\Q)-E(\Z[\frac{1}{R_D}]) : \hat{h}(P) \leq x \}.$$
 We find for large $x$, 
$$ |A_x| \approx \sum_{{{P \in E(\Q)}\atop{\hat{h}(P) \leq x}}} \hat{h}(P)^{-\delta}. $$
We now pick a basis $\{ P_i \}_{i=1}^r$ for the free part of $E(\Q)$ and write any $P \in E(\Q)$ as $\sum \lambda_i P_i+T$ with $\lambda_i \in \Z$ and $T \in E(\Q)_{{\rm tor}}$. Then $\hat{h}(P) \approx ||\lambda||^2 \cdot \log c$ for some constant $c$, an 
the above sum becomes 
$$ |A_x| \approx \sum_{{\lambda \in \Z^r-\{0\}}\atop{||\lambda||^2 \leq x}} ||\lambda||^{-2\delta}. $$
We group terms with $||\lambda||=m$ for a fixed integer $m$:
$$ |A_x| \approx \sum_{m=1}^{\sqrt{x}} m^{r-1} \cdot m^{-2\delta}. $$
We let $x \rightarrow \infty$, and find that $A$ is finite if this sum converges, which happens exactly for $2\delta-r+1 > 1$, i.e., $\delta>r/2$. This can be attained for $N$ sufficiently large.  \qed

\vv

With $r=1$, this implies  that $B_n$ doesn't have a divisor in $R_D$ only for finitely many $n$ as soon as $|D| \geq 2$. Applying it to the isogenous curve $E'$, it implies the same for $\{A_\ast\}$ and hence $\{C_\ast\}$.

Actually, the primitive part of $B_n$ is of size at least $\hat{h}(P)^{0.6 \cdot n^2}$ (cf.\ Silverman \cite{SilvWief}, Lemma 9 for an estimate $\hat{h}(P)^{n^2/3}$ and \cite{Streng} for a proof with a factor 0.6, using elliptic transcendence theory). We can apply the same argument to the primitive part of $B_n$. Furthermore, taking the $ABC$-conjecture for granted, if $E$ has $j$-invariant $0$ or $1728$, then we even know that the squarefree primitive part of $B_n$ is of the same order (Lemma 13 in loc.\ cit.), and this gives a heuristical proof of $R_D$-odd-primitivity for $|D| \geq 2$ on such curves. 

\vv

(\ref{heur}.3) One might note the following about the error term in (\ref{heur}.1): Shanks \cite{Shanks} analysed (\ref{heur}.1) in case  $M$ is the complement of the set of sums of two squares (cf.\ Ramanujan's first letter to Hardy) and noted that the first two terms give an accuracy of only $0.005$ at $x=10^7$. 

\vv

(\ref{heur}.4) Note further that for $E$ having a rational 2-torsion point, $B_n$ can be prime only finitely often, as follows from \cite{Everest} (since it arises as image sequence under an isogeny).  It is actually conjectured (see loc.\ cit.) that $B_n$ can only be prime for $n \leq K$ and some constant $K$ independent of $E$ and $P$; this is related to the elliptic Lehmer problem. It is reasonable to expect that $B_n$ has $m$ distinct odd order primitive prime factors as soon as $n \geq K$ for some constant $K$ only depending on $P$ and $E$ and $m$ (and maybe even only $m$). Granting that the (many) prime factors of $B_n$ are equidistributed over residue classes, the probability that at least one of the them is inert in a given $\Q(\sqrt{d})$ is very high.

\vv

\paragraph {\bf Further remarks.} \label{remconj} (i) One can wonder whether the property of being $R_D$-primitive is very sensitive to the choice of $D$, so ask whether it is true that {\sl for every elliptic curve $E$ (respectively, such that $(0,0) \in E[2]$) and every non-empty set $D$ of discriminants, for every rational point $P$ of infinite order and sufficiently large height, $B_n$ (respectively $C_n$) satisfies the $R_D$-primitivity condition}. There is some evidence that the $R_D$-primitive part doesn't behave the same for all $D$. For example, if $E$ has complex multiplication by some $d 
\in D$, there appear to be ``more'' split primes. This is explained by a Zsigmondy's theorem for an interpolation of the usual elliptic divisibility sequence by a sequence indexed by all endomorphisms of the curve, see Streng \cite{Streng}. Another example is Rubin's proof in \ref{exa}.  

\vv

(ii)  It is interesting to observe that the multiplicative group (disguised as the Fibonacci sequence) played an essential r\^{o}le in
the original proof of ${\rm HTP}(\Z)$. 
However, the analogue of our conjectures for linear recurrent sequences or the multiplicative group, i.e., an ``inertial classical Zsigmondy's theorem'', is almost certainly false. Let us reason heuristically for the sequence $\{a^n-1\}_{n \geq 1}$ for fixed $a$. The probability that $a^n-1$ is divisible only by primes outside $R_D$ is $[\log(a^n-1)]^{-\delta}$ with $\delta=1-1/2^{|D|}$ (cf.\ (\ref{heur}.1)), so the number of $n \leq x$ for which this holds is 
approximately $$\sum\limits_{n \leq x} \log[(a^n-1)]^{-\delta} \approx \ \sum\limits_{n \leq x} n^{-\delta}$$ which diverges if $x \rightarrow \infty$ for all $\delta$. Also, 
a general term of such a sequence (if $a$ is not composite) can be prime infinitely often. This is why we really need elliptic curves.

\vv

(iii) It is easy to formulate an analogue of the above conjecture for elliptic curves over global function fields. Especially in the case of an isotrivial
curve (e.g., the ``Manin-Denef curve'' $f(t)y^2=f(x)$), some information 
can be found in the literature, cf.\ \cite{Pheidas}.

Another function field analogue of \ref{IZ} is the following: let $\phi$ be a rank-2 $\F_q[T]$-Drinfeld module over $\F_q(T)$ (see, e.g.\ \cite{Goss}). If $x \in \F_q[T]$ is a polynomial of sufficiently large degree with $\phi_a(x) \neq 0$ for all $a \in \F_q[T]$, then for all polynomials $n$, $\phi_n(x)$ is divisible by an irreducible polynomial $\wp$ coprime to $\phi_m(x)$ for all $m$ of degree $\deg(m)<\deg(n)$, such that $\wp$ is inert in at least one of $\F(T)(\sqrt{d})$ for $d$ in a finite set of polynomials. A Drinfeld module analogue of Zsigmondy's theorem was proven by Hsia (\cite{Hsia}). 

\vv

(iv) The weaker statement that every $C_n/C_1$ has an odd order divisor from $R_d$, but not necessarily primitive, is equivalent to the fact that each of the ``fibrations in conics over $E$'' $$C/C_1=f(X,Y) \wedge C^2=A^4+aA^2B^2+bB^4$$ has only finitely many rational (${\bf P}^1$-)fibres over $E$, where $f$ runs over the classes of binary quadratic forms of the correct discriminant. For example, since $\Q(\sqrt{5})$ has class number one, related to example \ref{exa} is the diophantine equation
$$ (A^2+B^2)(A^2+11B^2)= 3^2 \cdot 5^2 \cdot (X^2-5Y^2)^2, $$
a smooth projective K3-surface whose  rational points should be found. 
One is reminded of the trouble deciding whether or not Martin Davis's equation has finitely many solutions (cf.\ Shanks and Wagstaff \cite{Shanks2}, again using Landau-Serre type estimates).

\vv

(v) In conjecture \ref{IZ}, one can move the point $P$ to $(0,0)$ by a rational change of coordinates. Then the conjecture becomes purely a statement about the division points on $E$, as we then have $$B_n^2 = \pm n^2 \prod_{Q \in E[n]} x(Q).$$ 
\vv

(vi) If $E$ has complex multiplication, then one has a divisibility sequence $\{ B_\alpha \}$ associated to any $\alpha \in \mbox{End}(E)$ (cf.\ \cite{Chud}). A similar theory with similar conjectures can be worked out. For the analogue of Zsigmondy's theorem, see \cite{Streng}.

\vv

\vv

\sectioning{A density version of the conjecture.} \label{densversion}

\vv

\paragraph {\bf Periodicity technique.} \ There is a principle of periodicity of elliptic divisibility sequences that can be used to prove density versions of the conjectures. Here is an example for Conjecture \ref{IZ}: the point $P=(-2,4)$ is non-singular modulo all primes and of infinite order on the curve $E \ : \ y^2=x^3+7x^2+2x$. The sequence $\{B_\ast\}$ for $P$ starts as $(1,2^2,3\cdot 11,2^3 \cdot 5^2, \dots)$ up to signs. 

\vv

\paragraph {\bf Definition.} \ The {\sl rank of apparition} $\rho_p=\rho_p(X_\ast)$ of a prime $p$ in a sequence $\{ X_\ast \}$ is the smallest $n$ for which $p|X_n$. 

\vv

\paragraph {\bf Periodicity} \label{Per} (Morgan Ward \cite{Ward:48}, section III). {\sl Assume that the sign of $B_n$ is chosen so that $B_n=\psi_n(P)$ for the classical division polynomial $\psi_n$. Assume $p>3$ has rank of apparition $\rho_p>3$ in $\{ B_\ast \}$. Then that sequence is periodic with period $\pi_p$ given by 
$$ \pi_p = \rho_p \cdot 2^{a_p} \cdot \tau_p, $$
where $\tau_p$ is the least common multiple of the (multiplicative) orders $\epsilon$ and $\kappa$ of $B_{\rho_p-1}$ and $B_{\rho_p-2}/B_2$ modulo $p$, respectively; and where $a_p=1$ if both $\epsilon$ and $\kappa$ are odd, $a_p=-1$ if both $\epsilon$ and $\kappa$ are divisible by the same power of $2$, and $a_p=0$ otherwise.}

\vv

We now look at the behaviour of the Jacobi symbol of $B_\ast$ modulo a given prime $p$. To avoid sign problems, we choose $p = 1$ mod $4$. For example, our sequence is periodic
modulo $5$ with period 8. Hence the sequence of Jacobi symbols $(\frac{5}{B_n})=(\frac{B_n}{5})$ (by quadratic reciprocity) is periodic with the same period, and its repeats $$(1,1,-1,0,-1,1,1,0) \mbox{ mod } 5.$$ This implies that 
$(\frac{5}{B_n})=-1$ whenever $n=\pm 3$ mod $8$, so all $B_s$ for $s$ a prime congruent to $\pm 3$ mod $8$ have a primitive odd order divisor in $R_{5}$.

In this case, one can do a little better. Assume $n=s^e$ is a power of a prime $s= \pm 3$ mod $8$. Then for $e$ even, $s^e = 1$ mod $8$ and $s^{e-1} = \pm 3$ mod $8$, whereas for $e$ odd, we have $s^e = \pm 3$ mod $8$ and $s^{e-1} = 1$ mod $8$. From periodicity, we see that in any case the Jacobi symbol of $B_{s^e}/B_{s^{e-1}}$ is $-1$, so the number is divisible by an odd order divisor in $R_{5}$. We conclude:

\vv

\paragraph {\bf Proposition.} \ \label{exadense} {\sl If $\{ B_\ast \}$ is the elliptic divisibility sequence associated to $(2,-4)$ on $y^2=x^3+7x^2+2x$, then any $B_{s^e}$ for $s$ a prime number $=\pm 3$ mod 8 has a {\sl primitive} odd order divisor from $R_5$. In particular, the set $\{ s \mbox{ prime } : \ B_s \mbox{ has a primitive odd order divisor from } R_5 \}$ has Dirichlet density at least $2/\varphi(8)=1/2$.} \qed

\vv

One can go on and create a race between inertial conditions in different $\Q(\sqrt{p})$ and the period of $\{B_\ast\}$ modulo $p$. We do this for the first few $p=$ 1 mod $5$ and the above curve and point (leaving out the easy computations). For $p=13$, the sequence has period 36 and for $s = \pm 5,7,11,13$ mod $36$, $(\frac{B_s}{13})=-1$. For $p=17$, all Kronecker symbols are positive. For $p=29$, the period is 38, and $s= \pm 9,11,15$ mod $38$ give a negative Kronecker symbol. For $p=37$, no new residue classes occur. For $p=41$, the period is $42$, and $s=\pm 13, 17$ mod $42$ have negative Kronecker symbol. For $53$, the period is $66$ and $s=\pm 5,7,25,29$ mod $66$ have negative Kronecker symbol. An easy density computation gives:

\vv

\paragraph {\bf Proposition.} \ \label{exadense2} {\sl Let  $\{ B_\ast \}$ denote the elliptic divisibility sequence associated to $(2,-4)$ on $y^2=x^3+7x^2+2x$, and let $D=\{5,13,29,41,53\}$. Then the set $$\{ s \mbox{ prime}\; :\; B_s \mbox{ has a primitive odd order divisor from } R_D \}$$ has Dirichlet density at least $43/45 \geq 95.5 \%$.} \qed

\vv

\vv

\paragraph{{\bf Remark.}} \ It is an interesting question whether, given any elliptic divisibility sequence $B_\ast$, and $\varepsilon>0$, one can choose a set $D$ such that $$\{ s \mbox{ prime}\; :\; B_s \mbox{ has a primitive odd order divisor from } R_D \}$$ has density at least $1-\varepsilon$.




\def\refname{\noindent\normalsize{References}}

\begin{footnotesize} 

\bibliographystyle{plain}

\vv

\vv

\noindent Mathematisch Instituut,
Universiteit Utrecht, Postbus 80010, 3508 TA  Utrecht, Nederland 

\noindent Email: {\tt cornelissen@math.uu.nl}

\vv

\noindent Departement Wiskunde, Statistiek en Actuariaat, Universiteit Antwerpen, Prinsstraat 13, 2000 Antwerpen, Belgi{\"e}

\noindent Email: {\tt zahidi@logique.jussieu.fr}

\end{footnotesize}

\end{document}